\documentclass[12pt,a4paper]{article}

\usepackage{amssymb}
\usepackage{amsmath}
\usepackage{subfigure}
\usepackage{color}
\numberwithin{equation}{section}
\usepackage[figuresright]{rotating}

\newtheorem{theorem}{Theorem}[section] 
 
\newtheorem{defi}[theorem]{Definition}
\newtheorem{rem}[theorem]{Remark}
\newtheorem{ass}[theorem]{Assumption}

\newcommand{\setR}{\mathbb{R}} 
\newcommand{\setC}{\mathbb{C}}

\newcommand{\setN}{\mathbb{N}}

\newcommand{\bu}{\mathbf{u}}
\newcommand{\bv}{\mathbf{v}}
\newcommand{\bw}{\mathbf{w}}

\newcommand{\bV}{\mathbf{V}}
\newcommand{\bfc}{\mathbf{f}}
\newcommand{\bg}{\mathbf{g}}
\newcommand{\bsigma}{\mathbf{\sigma}}
\newcommand{\bepsilon}{\mathbf{\epsilon}}
\newcommand{\bn}{\mathbf{n}}
\newcommand{\bx}{\mathbf{x}}
\newcommand{\by}{\mathbf{y}}
\newcommand{\smat}[1]{\left(\begin{smallmatrix} #1 \end{smallmatrix}\right)}

\DeclareMathOperator{\dv}{div}

\DeclareMathOperator{\Bdv}{\mathcal{B}}
\DeclareMathOperator{\sign}{sign}

\DeclareMathOperator{\spann}{span}
\newcommand{\bui}{\bu^{\rm int}}
\newcommand{\bvi}{\bv^{\rm int}}
\newcommand{\bue}{\bu^{\rm ext}}
\newcommand{\bve}{\bv^{\rm ext}}

\newcommand{\mpara}{\zeta}
\newcommand{\WGx}{\xi}
\newcommand{\WGy}{\eta}
\DeclareMathOperator{\LT}{{\cal L}}
\newcommand{\iu}{\mathfrak{i}}	
\newcommand{\pol}{s_0}	
\newcommand{\poli}{s_1}	
\newcommand{\setWN}{\mathbb{S}}
\newcommand{\Oi}{\Omega_{\rm int}} 
\newcommand{\Oe}{\Omega_{\rm ext}}
\newcommand{\WGTH}{R}
\newcommand{\surf}{\Upsilon}
\newcommand{\WG}{W}
\newcommand{\nrWG}{L}
\newcommand{\emotion}{\mathcal{E}}
\newcommand{\Boundary}{{\partial \Omega}}
\newcommand{\BD}{\Boundary_D}
\newcommand{\BN}{\Boundary_N}
\newcommand{\curve}{\Gamma}
\newcommand{\Ba}[2]{a\!\left(#1,#2\right)}
\newcommand{\Bb}[2]{b\!\left(#1,#2\right)}

\newcommand{\Bae}[2]{a_{\rm ext}\!\left(#1,#2\right)}
\newcommand{\Bbe}[2]{b_{\rm ext}\!\left(#1,#2\right)}
\newcommand{\Baell}[2]{a^{\rm ext}_\ell\!\left(#1,#2\right)}
\newcommand{\Bbell}[2]{b^{\rm ext}_\ell\!\left(#1,#2\right)}
\newcommand{\Bai}[2]{a_{\rm int}\left(#1,#2\right)}
\newcommand{\Bbi}[2]{b_{\rm int}\left(#1,#2\right)}
\newcommand{\Baux}[2]{f_{\rm a}\left(#1,#2\right)}
\newcommand{\Lf}[1]{f\!\left(#1\right)}
\newcommand{\bai}{\varphi^{\rm int}}
\newcommand{\bay}{\varphi^{\rm trans}}
\newcommand{\bax}{\varphi^{\rm long}}
\newcommand{\bayell}{\varphi^{\rm trans,\ell}}

\newcommand{\Ma}{A}
\newcommand{\Mb}{B}
\newcommand{\Mai}{A^{\rm int}}
\newcommand{\Mbi}{B^{\rm int}}
\newcommand{\Mae}{A^{\rm ext}}
\newcommand{\Mbe}{B^{\rm ext}}
\newcommand{\My}{M_{\rm trans}}
\newcommand{\Dy}{D_{\rm trans}}
\newcommand{\Sy}{S_{\rm trans}}
\newcommand{\Mx}{M_{\rm long}}
\newcommand{\Dx}{D_{\rm long}}
\newcommand{\Sx}{S_{\rm long}}

\newcommand{\dfn}[1]{{#1}}
\newcommand{\dfo}[1]{\textcolor{blue}{}}

\title{Hardy space infinite elements for time-harmonic two-dimensional elastic waveguide
problems\footnote{Support from the Austrian Science Fund (FWF) under
grants W1245-N25 and P26252-N25 is acknowledged.}}
\author{Martin~Halla\footnote{Institute for Analysis and Scientific Computing,
Vienna University of Technology, Austria (\texttt{martin.halla@tuwien.ac.at})},
Lothar~Nannen\footnote{Institute for Analysis and Scientific Computing,
Vienna University of Technology, Austria (\texttt{lothar.nannen@tuwien.ac.at})}}

\begin{document}
\maketitle

\begin{abstract}
We consider time-harmonic linear elasticity equations in domains containing two-dimensional semi-infinite strips.
\dfo{Such problems exhibit solutions for which the signs of group and phase velocities differ.}
\dfo{tandard transparent boundary conditions, e.g. the}
\dfn{Since for such problems there exist modes with different signs of group and phase velocity, standard perfectly matched layer (PML) as well as standard Hardy space infinite element methods fail.} \dfo{perfectly matched layer (PML) method, select modes with positive phase velocity and hence yield unphysical solutions for such problems.}

We apply a recently developed infinite element method for a physically correct discretization of such waveguide problems which is based on a Laplace transform in propagation direction. In the Laplace domain the space of transformed
solutions can be separated into a sum of a space of incoming and a space of outgoing functions where both function spaces are certain Hardy spaces. The Hardy space is chosen such that the construction of a simple infinite element is
possible.

The method does not use a modal separation and works on intervals of frequencies. On those intervals the involved
operators are frequency independent and hence lead to linear eigenvalue problems when computing resonances.
Numerical experiments containing convergence tests and resonance problems are included.
\end{abstract}


\section{Introduction}\label{sec:intro}
Computational methods for wave equations bear great attention due to their huge importance for real live problems. In solid mechanics the application reaches from simulating seismic waves,
non-destructive testing to material characterization~\cite{Chimenti:97,Langenberg:12}.
Such equations are usually posed on unbounded domains $\Omega$. Mesh based methods, i.e. finite difference/volume/element methods deal with that difficulty by truncation of 
$\Omega$ to a bounded subdomain $\Oi$ and special treatment of the exterior domain $\Oe=\Omega\setminus\overline{\Oi}$. 

Based on representation formulas for the solution in $\Oe$ with given boundary data at the artificial boundary $\surf=\overline{\Oe}\cap \overline{\Oi}$ local approximations of the
Dirichlet-to-Neumann operator at $\surf$ can be used (see \cite{Givoli:04} for a review). Also based on
representation formulas for the solution (via a Green function) boundary element methods typically lead to non-local but more accurate approximations.

A method, which does not directly use a representation formula, is the complex scaling method, reintroduced by B\'erenger \cite{Berenger:94} as perfectly matched layer method (PML) 
for electromagnetic waves. 
The method became soon very popular and can be classified as todays standard method for treating unbounded domains. In~\cite{HohageNannen:09} a Hardy space infinite element method (HSIE) 
was introduced for Helmholtz problems. The method can be understood as a special infinite element method, relying on the pole condition~\cite{PC1}. It features many advantages of PMLs, 
in particular it does not depend directly on a representation formula. However, the theoretical background of this method is different to that of the PML.

As for all linear wave equations, a fundamental technique in understanding linear elastic wave equations is the study of most simple solutions~\cite{Graff,Achenbach73}. In general unbounded domains 
these are plane waves $e^{\iu\mathbf{\kappa}\cdot\bx -\iu\omega t}\bw$, where as in waveguides, such as plates $\setR^2\times I, I\subset\setR$ and cylinders $\setR\times D, D\subset \setR^2$, 
boundary conditions have to be respected by the solutions.
The general form in these cases are modal waves. In semi-infinite cylinders they take the form $e^{\iu\kappa x-\iu\omega t}\bw(\by)$, 
where $t\in\setR_{>0}$ is the time variable, $x\in\setR_{>0}$ is the \dfn{longitudinal}\dfo{unbounded} coordinate, $\by\in D$ the vector of \dfn{transverse}\dfo{bounded} coordinates, $\omega\in\setR_{>0}$ the angular frequency, $\kappa$ the wavenumber 
and $\bw/|\bw|$ the direction of displacement. The frequency $\omega$ and the wavenumber $\kappa$ have to fulfill a wave equation specific 
dispersion relation $\omega(\mathbf{\kappa})$ to yield a solution. The wave travels in direction $\sign(\kappa)$ with phase velocity $\omega/\kappa$, 
whereas the energy is transported in direction $\sign(\partial\omega(\kappa)/\partial\kappa)$ with group velocity $\partial\omega(\kappa)/\partial\kappa$~\cite{Lighthill:65}.

The possibility of waves with different signs of group and phase velocity was already discussed by Lamb in 1904~\cite{Lamb:1904}. 
Such mismatches between the directions of propagation of phase and energy cannot only happen in waveguides, but also in anisotropic materials~\cite{Joly:03}. 
Although PML methods select waves by their propagation direction of phase, the attention to mismatches between propagation directions of phase and energy in elastic materials 
was omitted for long~\cite{Joly:03}. This may result in a wrong selection of outgoing/incoming waves and can lead to an exponential growth of solutions in the damping layer of PMLs. 
Such cases are usually reported as instabilities, due to the definition of stability, but first of all, in such cases an unphysical radiation condition is incorporated in the PML formulation.

Sometimes, such problems can be solved using a transformation of variables~\cite{BNiBonLeg:2004,Abarbanel:99}. 
In~\cite{SMARTlayer1,SMARTlayer2} SMART layers were introduced for seismic waves. 
They allow more flexibility in the damping \dfn{than}\dfo{then} PMLs, but lose the property of being perfectly matched.
For semi-infinite linear elastic cylinders in the time-harmonic setting, two methods based on bi-orthogonal relations of modal solutions were proposed: 
The method presented in~\cite{Baronian} is based directly on a modal representation of the solution,  whereas in~\cite{Skeltonetal:07, Bonnet} ways to modify standard PMLs are reported.
They rely on a smart post processing by exchanging backward incoming with backward outgoing modes.
\dfo{However, these methods depend non-linearly on $\omega^2$.}
\dfn{Since, these methods involve modal solutions, which depend non-linearly on $\omega^2$, they depend
themselves non-linearly on $\omega^2$.}
\dfn{Thus if they are used to compute resonances, they lead to non-linear eigenvalue problems. Although solvers for non-linear eigenvalue problems exist (e.g.~\cite{NEVP}),
they are more involved than solvers for linear eigenvalue problems.}

In~\cite{HallaHohageNannenSchoeberl:14} a new family of Hardy space infinite elements was introduced. It is based on a generalized pole condition, which allows for different signs of phase velocities.
The achievement of this paper is to pick up these results and apply it to two-dimensional time-harmonic semi-infinite cylindrical linear elastic waveguide problems. 
A big advantage of this method is, that \dfn{$\omega^2$}\dfo{the square of the frequency} only enters as a scalar coefficient in the method.
\dfn{Hence, the discretization of a resonance problem leads to a generalized linear matrix eigenvalue problem,
which can be treated with a standard solver.}
\dfo{This leads to linear eigenvalue problems when looking for solutions
to the resonance problem.}
Moreover, the numerical results indicate a super-algebraic convergence with respect to the number of degrees of freedom in the \dfn{longitudinal}\dfo{infinite} direction.

The outline of the paper is as follows: Section~\ref{sec:setting} formulates the diffraction problem to solve, in particular the modal radiation condition
and its reformulation as \dfn{a} pole condition.
\dfn{In Section~\ref{sec:HSM} we introduce the Hardy space infinite element in one dimension as in \cite{HallaHohageNannenSchoeberl:14} and discuss the choice of method parameters for the investigated elasticity
problem. Section~\ref{sec:disc-method} explains the discretization of the elasticity problem with the use of
tensor product basis functions.}
\dfo{In Section~\ref{sec:disc-method} we describe a (rather general) infinite element discretization.
In Section~\ref{sec:HSM} we specify the Hardy space infinite element, state its element matrices and discuss the choice of method parameters.}
Section~\ref{sec:resprob} deals with the \dfn{spectral objects and properties}
\dfo{definition and interpretation} of resonance problems.
Finally, in Section~\ref{sec:numres} we give numerical examples. 



\section{General setting}\label{sec:setting}

\subsection{Geometry and elasticity Equations}
\begin{figure}[tb]
\begin{center}
\resizebox{0.6\textwidth}{!}{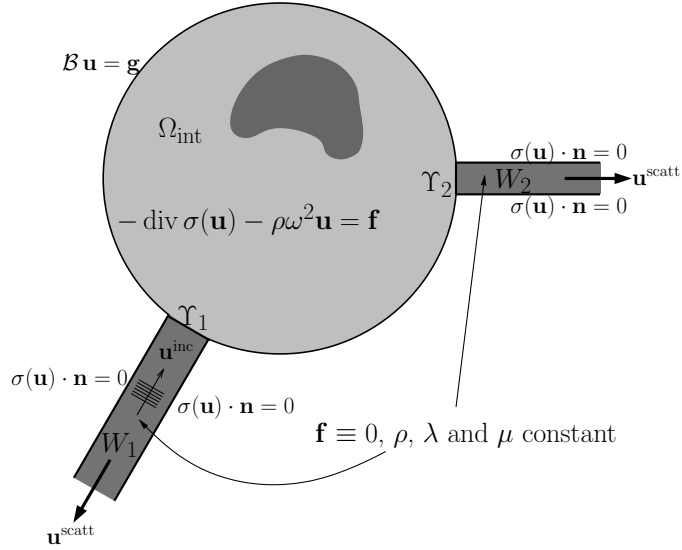}
\end{center}
\caption{\dfn{sketch of an elastic waveguide problem under consideration in this paper}}
\label{Fig:WGproblem}
\end{figure}

\dfn{We are considering in this paper elastic waveguide problems in two dimensions (see Fig.~\ref{Fig:WGproblem} for a typical situation): The domain of interest $\Omega$ consists of a bounded interior domain $\Oi$ and $\nrWG$ semi-infinite waveguides $\WG_1,\dots,\WG_\nrWG$ with interfaces $\surf_1,\dots,\surf_\nrWG$.}
\dfo{Let $\Omega=\Oi\cup\bigcup_{\ell=1}^\nrWG(\WG_\ell\cup\surf_\ell)$ be a Lipschitz domain, which is the disjoint union of
a bounded Lipschitz domain $\Oi\subset\setR^2$, $\nrWG$ semi-infinite strips (waveguides) $\WG_1,\dots,\WG_\nrWG$ and
interfaces $\surf_1,\dots,\surf_\nrWG$. More precisely, the waveguides $\WG_\ell$ and interfaces $\surf_\ell$ are of the form
$\WG_\ell=\emotion_\ell\left((0,\infty)\times\tilde\surf_\ell\right)$ and
$\surf_\ell=\emotion_\ell\left(\{0\}\times\tilde\surf_\ell\right)$ where $\emotion_\ell\colon\setR^2\to\setR^2$ is an Euclidean motion and $\tilde\surf_\ell:=(-\WGTH_\ell,\WGTH_\ell)$ with $\WGTH_\ell>0$. }
The two-dimensional time-harmonic isotropic linear elasticity problem is given by
\begin{subequations}\label{eq:lame}
\begin{align}\label{eq:lame-inner}
-\dv \bsigma (\bu) -\rho\omega^2 \bu &= \bfc \text{ in } \Omega,\\
\Bdv \bu &=\bg \text{ on } \partial \Omega,\\
\bu-\bu^{\rm inc} &\text{ satisfies a radiation condition in each waveguide } \WG_\ell,\quad \ell=1,\dots,\nrWG. \label{eq:abstrRC}
\end{align}
\end{subequations}
Here, $\Re (\bu(\bx)e^{-\iu\omega t})$ for $\bx\in\Omega$ and time $t>0$ is the time-harmonic displacement vector,
$\rho>0$ the density, $\omega>0$ the angular frequency, $\lambda,\mu>0$ are Lam\'e parameters, $\bepsilon (\bu)=\frac{1}{2}(\nabla \bu + (\nabla \bu)^\top)$ 
the strain tensor, $\bsigma (\bu)=\lambda\dv \bu \cdot \mathbf{Id} +2\mu\bepsilon(\bu)$ the stress tensor, and $\bfc$ a
volumetric force\dfo{ with compact support in $\Oi$}. \dfn{In this paper we will always use the plain strain model (see e.g. \cite{Braess:FE}), i.e. the Lam\'e parameters are related to Young's modulus $E$ and Poisson's ratio $\nu$ by $\mu=\frac{E}{2(1+\nu)}$ and $\lambda=\frac{E\nu}{(1+\nu)(1-2\nu)}$.
Boundary conditions at $\partial\Oi\cap\partial\Omega$ are formulated in terms of a trace operator $\Bdv$,  e.g. a Dirichlet trace operator $\Bdv \bu=\bu$ or a Neumann trace operator $\Bdv \bu=\bsigma(\bu)\cdot \bn$ with outer normal vector $\bn$, and a boundary datum $\bg$.}

We will assume in the following, that $\rho$, $\lambda$ and $\mu$ are constant in each waveguide $\WG_\ell$ \dfn{, that there exists no volumetric force $\bfc$ in the waveguides, and that we have traction free boundary conditions $\bsigma(\bu)\cdot \bn = 0$ at the boundaries $\partial \WG_\ell \cap \partial \Omega$.}
\dfo{$\Bdv$ denotes a trace operator, e.g. for $\partial\Omega=\BD\cup\BN$ the Dirichlet trace operator $\Bdv_{\BD} \bu=\bu$
on $\BD$ and the Neumann trace operator $\Bdv_{\BN} \bu=\bsigma(\bu)\cdot \bn$ on $\BN$ with outer normal vector $\bn$
on $\partial\Omega$. $\bg$ is a boundary datum with compact support in $\partial\Oi\cap\partial\Omega$.} The terms
radiation condition and incident field $\bu^{\rm inc}$ in \eqref{eq:abstrRC} \dfn{and in Fig.~\ref{Fig:WGproblem}} will \dfo{ be defined later on in
Def.~\ref{defi:modal_rad_cond} and in Def.~\ref{def:PoleCondition}.} \dfn{become clear in the next subsection.}

\dfo{We will also consider resonance problems, which have the form \eqref{eq:lame}, but $\bfc ,\bg$, and $\bu^{\rm inc}$
vanish, $\omega$ may be complex valued with positive real part, and both $\omega$ and $\bu\neq 0$ are considered as
unknowns.}
\dfn{In addition to the scattering problem \eqref{eq:lame}, where the angular frequency $\omega$ and the sources $\bfc ,\bg$, and $\bu^{\rm inc}$ are given, we will consider the corresponding resonance problem:
Find resonances $\omega\in \setC$ with positive real part and non-trivial resonance functions $\bu$
such that
\begin{subequations}\label{eq:res_classic}
\begin{align}
-\dv \bsigma (\bu) -\rho\omega^2 \bu &= \mathbf{0} \text{ in } \Omega,\\
\Bdv \bu &=\mathbf{0} \text{ on } \partial \Omega,\\
\label{eq:RC_resprob} \bu  &\text{ satisfies a radiation condition in each waveguide } \WG_\ell,\quad \ell=1,\dots,\nrWG.
\end{align}
\end{subequations}}
\dfo{In the case of the plain strain model the Lam\'e parameters are connected with Young's modulus $E$ and Poisson's ratio $\nu$ through the relations
\begin{subequations}
\begin{equation}
\mu=\frac{E}{2(1+\nu)}, \quad \lambda=\frac{E\nu}{(1+\nu)(1-2\nu)}.
\end{equation}
In the case of the plain stress model the relations are
\begin{equation}
\mu=\frac{E}{2(1+\nu)}, \quad \lambda=\frac{E\nu}{(1+\nu)(1-\nu)}.
\end{equation}
\end{subequations} 
See \cite{Braess:FE} for further details. In this paper $\mu, \lambda$ will always be related to $E, \nu$ through the plain strain model.}

\subsection{Radiation condition}\label{Sec:ModalAnalysis}
\begin{figure}
\centering
\subfigure[\label{fig:disprel_martin2}]{\resizebox{0.45\textwidth}{!}{\includegraphics{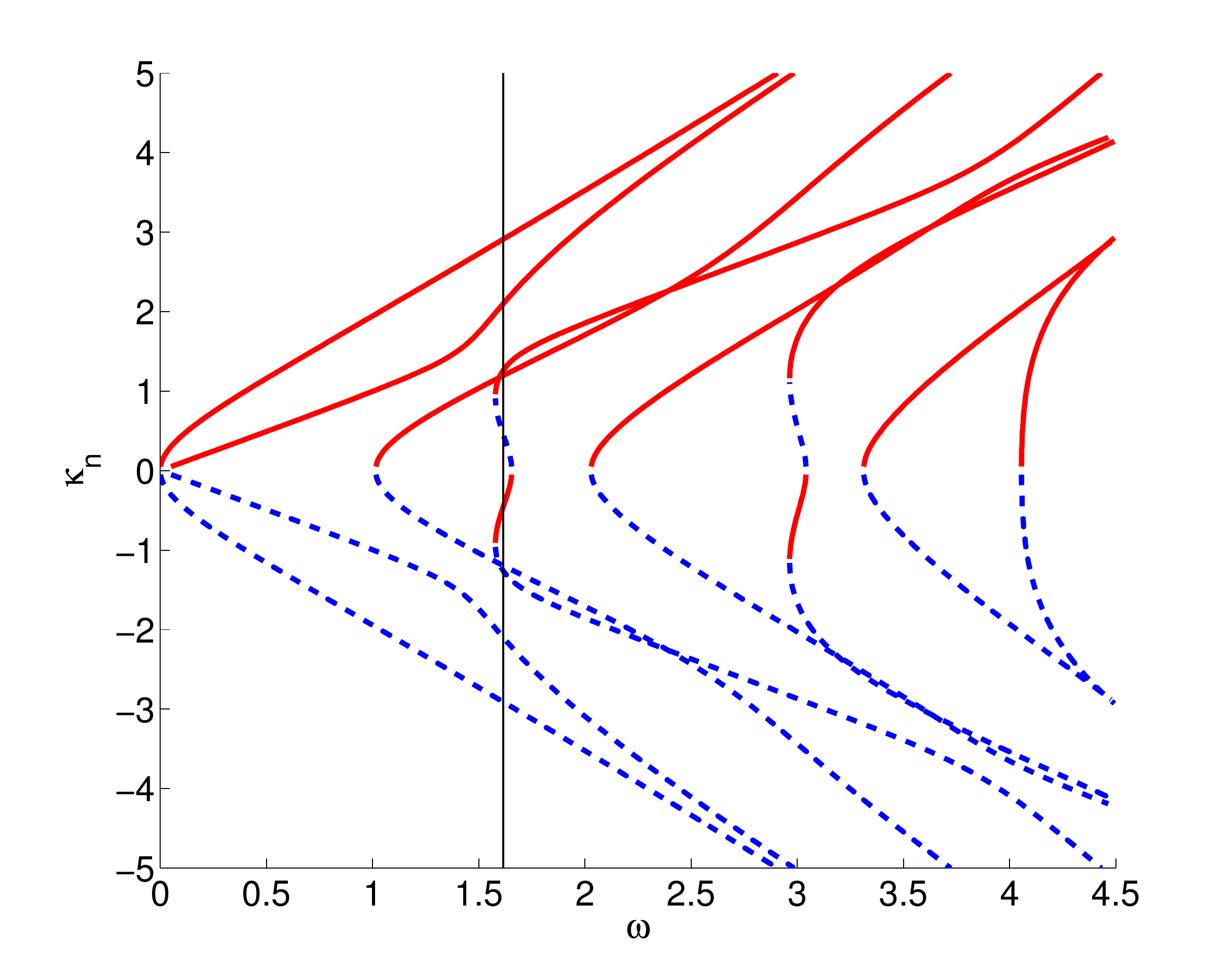}}} \hfill
\subfigure[\label{fig:wn3}]{\resizebox{0.45\textwidth}{!}{\includegraphics{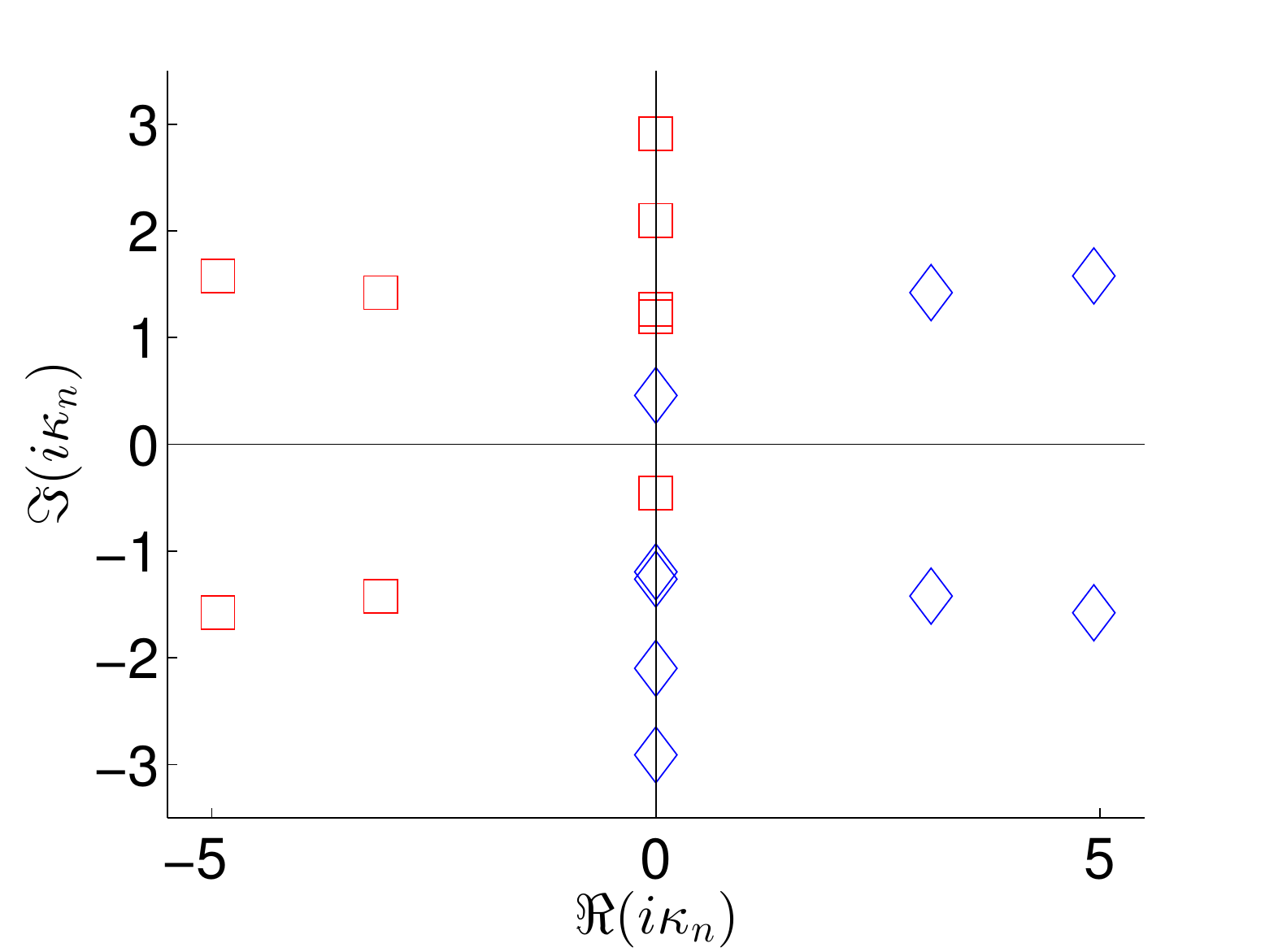}}} 
\caption{left: the first eight dispersion curves. Modes corresponding to the red solid part have positive group velocity and modes corresponding to the blue dashed part have negative group velocity.
right: the first nine outgoing/incoming wavenumbers multiplied with $\iu$ marked with red squares/blue diamonds at frequency $\omega=1.615$. Parameters are $H=\rho=E=1$, and $\nu=0.2$.}
\label{fig:wn}
\end{figure}

In order to define a physically correct radiation condition \eqref{eq:abstrRC}, we consider in this section one single
waveguide. The generalization to multiple waveguides is straightforward. Let $\WG:=\setR^+\times (-\WGTH,\WGTH)$ be a
reference waveguide and $\rho$, $\lambda$ and $\mu$ be constant and positive. As for most other problems the radiation
condition can be derived by an analytic representation of solutions. For waveguides a convenient representation is a
modal sum. For \dfn{given angular frequency} $\omega>0$ we call $\bu(\bullet,\bullet;\omega)$ a mode with wavenumber $\kappa(\omega)\in \setC$, if it
has the form
\begin{align}\label{eq:mode}
 \bu(\WGx,\WGy;\omega)=e^{\iu \kappa(\omega) \WGx} \bw(\WGy;\omega), \qquad (\WGx,\WGy) \in \WG, \omega\in\setR^+,
\end{align}
and solves
\begin{subequations}\label{eq:modal}
\begin{align}
-\dv \bsigma (\bu )-\rho\omega^2 \bu &= 0,\qquad &&(\WGx,\WGy) \in \WG, \label{eq:modal_de}\\
\bsigma(\bu)\cdot\smat{0\\1}  &=0,\qquad &&(\WGx,\WGy) \in \setR_+\times \{-\WGTH,\WGTH\}.\label{eq:modal_bcWG}
\end{align}
\end{subequations}
It is straightforward to see that if $e^{\iu \kappa(\omega) \WGx} \smat{\bw_1(\WGy;\omega)\\ \bw_2(\WGy;\omega)}$ is a
mode, so is $e^{-\iu \kappa(\omega) \WGx} \smat{\bw_1(\WGy;\omega)\\-\bw_2(\WGy;\omega)}$ as well as
$e^{-\iu\overline{\kappa}(\omega) \WGx} \smat{\overline{\bw}_1(\WGy;\omega)\\\overline{\bw}_2(\WGy;\omega)}$ and $e^{\iu \overline{\kappa}(\omega) \WGx} \smat{\overline{\bw}_1(\WGy;\omega)\\-\overline{\bw}_2(\WGy;\omega)}$. 
A physical solution should be bounded for $\WGx \to \infty$. Hence, if $\Im(\kappa(\omega))\neq0$ we want to exclude
modes with $\Im(\kappa(\omega))<0$. In this case, we call the exponentially decaying modes \emph{evanescent}. If
$\Im(\kappa(\omega))=0$, both modes with wavenumbers $\pm\kappa(\omega)$ stay bounded and it is not obvious which one is
physically relevant. We give here two approaches leading to the same distinguishing criterion. For both we assume, that $\Re(\partial_{\omega}\kappa(\omega))\neq0$, $\Im(\partial_{\omega}\kappa(\omega))=0$.
\begin{enumerate}
 \item \cite{Lighthill:65} states that the velocity of energy transport of a mode is \dfo{proportional to}\dfn{given by} the group velocity
 $\partial\omega(\kappa)/\partial \kappa$ which has the same sign as $\partial \kappa(\omega)/\partial \omega$. Since in
 problem~\eqref{eq:modal} there is no source in $\WG$, energy should be radiated to infinity and therefore
 $\partial\kappa(\omega)/\partial \omega>0$. See also \cite[Rem. 3.1]{Orazov}.
 \item For the \emph{limiting absorption principle} we add an artificial small damping $\epsilon>0$ to the system, i.e.
 we substitute $\omega>0$ by $\omega+\iu\epsilon$. The corresponding solution $\bu_\epsilon$ to the damped version of
 \eqref{eq:modal} should be bounded for $x \to \infty$. 
 Since by linearization $|\exp\left(\iu \kappa(\omega+\iu\epsilon)\xi\right)|\approx\exp\left(-\frac{\partial\kappa(\omega)}{\partial \omega}\xi\right)$ this leads to 
 $\partial\kappa(\omega)/\partial \omega>0$. For further details see \cite[Cor. 3.1]{Orazov}.
\end{enumerate}
Therefore, if $\partial\kappa(\omega)/\partial \omega>0$ we call 
$e^{\iu \kappa(\omega) \WGx} \smat{\bw_1(\WGy;\omega)\\ \bw_2(\WGy;\omega)}$ an \emph{outward propagating mode} and $e^{-\iu \kappa(\omega) \WGx} \smat{\bw_1(\WGy;\omega)\\ -\bw_2(\WGy;\omega)}$ 
an \emph{inward propagating mode}. 
Note, that there exist modes with different signs of phase and group velocity (see Fig.~\ref{fig:wn}). 
\dfn{We call a wavenumber $\kappa(\omega)$ \emph{outgoing}, if there exists an evanescent or outward propagating mode of the form \eqref{eq:mode} satisfying \eqref{eq:modal}. In other words $\kappa(\omega)$ is outgoing (\emph{incoming}), if
\begin{enumerate}
 \item $\partial_{\omega}\kappa(\omega)>0$ ($\partial_{\omega}\kappa(\omega)<0$) for real wavenumbers $\kappa(\omega)\in \setR$, and
 \item $\Im(\kappa(\omega))>0$ ($\Im(\kappa(\omega))<0$) for non-real wavenumbers $\kappa(\omega)\not\in \setR$.
\end{enumerate}
A function satisfies the radiation condition in the waveguide $\WG$ and is called \emph{outgoing}, if it can be approximated by a linear combination of evanescent and outward propagating modes.}
\dfn{It is not a priori clear, that the traces of evanescent and outward propagating modes on the waveguide interface $\surf=\{0\}\times (-R,R)$ are dense in $L^2(\surf)^2$. In the following Remark we cite some results from \cite{Orazov}, which justify our definition of the radiation condition.}
\begin{rem}\label{Rem:TheoModes}
In~\cite{Orazov} a rigorously analysis of modal decompositions for \dfn{two and} three dimensional elastic waveguides is given using a quadratic eigenvalue problem.
There, not only modes (eigenfunctions) of the form~\eqref{eq:mode} \dfn{are considered, but also
the generalized eigenspaces of the investigated quadratic eigenvalue problem (see \cite[(0.10)]{Orazov}).
We refer to these functions \dfn{of the generalized eigenspace, which are not eigenfunctions,} as associated modes in this paper}.
We briefly summarize the main results of~\cite{Orazov}: 
\begin{enumerate}
 \item \cite[Thm. 1.5]{Orazov}: For fixed $\omega>0$ there exist wavenumbers $\kappa_n(\omega),n\in\setN$, symmetrically situated relative to the real axis and the origin. They are situated in arbitrarily small angles, adjoining the imaginary axis, with the exception of a finite number of wavenumbers. In particular, only a finite number $\kappa_n(\omega)\in\setR$ exists.
 \item \cite[Thm. 3.8]{Orazov}: There exists a sequence of frequencies $0=\omega_1^2<\omega_2^2<\dots<\omega_n^2\to\infty$, such that for all $\omega\in\setR^+\setminus\{\omega_n,n\in\setN\}$ the implication $\kappa_n(\omega)\in\setR \Rightarrow \partial_{\omega} \kappa_n(\omega)\in\setR\setminus\{0\}$ holds.
 \item \cite[Thm. 2.5]{Orazov}: Assume $\omega\in\setR^+\setminus\{\omega_n,n\in\setN\}$. Then the \dfo{set}\dfn{traces} of modes and associated modes corresponding to $\{\kappa_n(\omega),n\in\setN\colon\Im(\kappa_n(\omega))>0\vee (\kappa_n(\omega)\in\setR\wedge\partial_{\omega}\kappa_n(\omega)>0)\}$ are dense and minimal in \dfo{$L^2(-\WGTH,\WGTH)$}\dfn{$L^2(\surf)^2$} as well as in \dfo{$H^1(-\WGTH,\WGTH)$}\dfn{$H^1(\surf)^2$}.
\end{enumerate}
\end{rem}
\dfo{
These results justify the following definition.
\begin{defi}[modal radiation condition]\label{defi:modal_rad_cond} 
Let $\omega\in\setR^+\setminus\{\omega_n,n\in\setN\}$ with $(\omega_n)_{n \in \setN}$ defined in Rem.~\ref{Rem:TheoModes}. We call a wavenumber $\kappa(\omega)$ \emph{outgoing} (or \emph{incoming}), if there exists a mode of the form \eqref{eq:mode} satisfying \eqref{eq:modal} and
\begin{enumerate}
 \item $\kappa(\omega)\in \setR$ and $\partial_{\omega}\kappa(\omega)>0$ (or $\partial_{\omega}\kappa(\omega)<0$) or
 \item $\Im(\kappa(\omega))>0$ (or $\Im(\kappa(\omega))<0$).
\end{enumerate}
We call a function $\bu\in [H^1_{\rm loc}(W)]^2$ \emph{outgoing}, if there exist $c_1,\dots,c_m\in\setC$, such that
$\|(\bu-\sum_{j=1}^m c_j \bv_j)\|_{\left(L^2(\{\xi\}\times(-R,R))\right)^2}\xrightarrow{\xi\to\infty}0$,
where the functions $\bv_j$ are the eigenfunctions and
associated functions to outgoing wavenumbers $\kappa_1(\omega),\dots,\kappa_n(\omega)\in\setR$. $\bu$ is called an
\emph{incident field} if it is a linear combination of eigenfunctions and associated functions to real incoming wavenumbers.
\end{defi}}

\subsection{\dfn{Numerical implementations of the radiation condition}}\label{sec:stateoftheart}
\dfn{In a numerical method for \eqref{eq:lame} the radiation condition as to be taken into account. Due to the simple waveguide geometry, methods based on modal decompositions as the one proposed in \cite{Baronian} are popular. Assuming that every solution in a waveguide can be expanded in a sum of modes $\sum_{n\in\setN}e^{\iu\kappa_n(\omega)\xi}\bw_n(\eta)$,
a Dirichlet-to-Neumann operator can be defined as
\begin{align}\label{eq:DtN}
\operatorname{DtN}\sum_{n\in\setN}e^{\iu\kappa_n(\omega)\xi}\bw_n(\eta)
:=\sum_{n\in\setN}\sigma_n\big(e^{\iu\kappa_n(\omega)\xi}\bw_n(\eta)\big).
\end{align}
Truncating the waveguides $W_l$ and posing
$\sigma_n(\bu)=\operatorname{DtN}\bu $
at the artificial ends leads to a formulation of~\eqref{eq:lame} in a bounded domain. For a numerical method standard finite element methods for the bounded domain can be used together with a discrete
Dirichlet-to-Neumann operator $(\operatorname{DtN})_N$, which can 
be constructed by using only a finite number of modes in~\eqref{eq:DtN}, i.e.
$(\operatorname{DtN})_N\sum_{n\leq N}e^{\iu\kappa_n(\omega)\xi}\bw_n(\eta)
:=\sum_{n\leq N}\sigma_n\big(e^{\iu\kappa_n(\omega)\xi}\bw_n(\eta)\big). $
It can be shown that the modes satisfy biorthogonal relations, which can be exploited
to implement $(\operatorname{DtN})_N$ in an elegant way, as done in \cite{Baronian}.
However, the method requires a precomputation of the modes. Moreover, since the modes depend on $\omega$, the discretization matrix will depend non-linearly on $\omega$. Hence, discretizing the resonance problem \eqref{eq:res_classic} leads to a large, non-linear eigenvalue problem. }

\dfn{In \cite{LM:08} a numerical method for a Helmholtz problem is proposed based on eigenfunction expansions of the bounded interior problem coupled to mode expansions in the waveguides. These approaches lead again
to non-linear, but comparatively small eigenvalue problems. Nevertheless, solving the non-linear eigenvalue problem numerically is a non-trivial task and of course a computation of the modes and interior eigenfunctions is needed.}

\dfn{A popular method which does not depend on the computation of the modes is the
perfectly matched layer method. The longitudinal direction $\WGx$ of the waveguide is thereby complex scaled by  $\alpha\in\setC$, i.e. $\tilde \WGx := \alpha \WGx$, leading to
a new equation for  $\tilde\bu(\WGx,\WGy):=\bu(\tilde \WGx,\WGy)$. If $\alpha$ is chosen such that $\Re(\alpha\iu\kappa_n(\omega))<0$ for all outgoing wavenumbers $\kappa_n(\omega)$
and $\bu$ is an outgoing solution of the original problem, than $\tilde\bu$ is exponentially decaying in $\WGx$. Hence, a
truncation of the infinite waveguide to a finite waveguide with zero Neumann or Dirichlet boundary condition at the
end introduces only a small error due to the exponential decay. The resulting equation can then be discretized
with standard methods, as it is posed on a finite domain.}
\dfn{The linear system matrix takes the form
$K_h-\omega^2M_h.$ 
When looking for resonances, the above leads to generalized linear matrix eigenvalue problems, which is a main advantage
in view of numerical algorithms.}

\dfn{The reason why perfectly matched layers are only of limited use for waveguide problems
with backward propagating modes, is that no $\alpha\in\setC$ can exist, such that $\Re(\alpha\iu\kappa_n(\omega))<0$
for all outgoing wavenumbers $\kappa_n(\omega)$ (see Fig.~\ref{fig:wn3}, a multiplication with $\alpha$ leads to a rotation of $\iu\kappa_n$ in the complex plane). The modified perfectly matched layer methods presented in \cite{Skeltonetal:07,Bonnet} combine complex scaling with a special treatment of the backward propagating mode. For scattering problems with only a few backward propagating modes this approach works
well, but again it leads to non-linear eigenvalue problems when discretizing the resonance problem \eqref{eq:res_classic}.
}

\subsection{Pole Condition}\label{Sec:PC}
Our goal \dfo{instead }is to construct a radiation condition,
\dfn{which combines the advantages of the methods presented in the previous subsection:
It should yield physically correct solutions and discretizations of resonance problems should
lead to linear matrix eigenvalue problems.}
\dfo{which is (at least for some domain of frequencies) independent of the frequency and the particular modal decomposition.}
Again for the sake of simplicity, we formulate the following only for the reference waveguide $\WG$. For $\omega>0$ let
$$\bu(\WGx,\WGy)=\sum_n^N c_n\bu_n(\WGx,\WGy), \qquad (\WGx,\WGy)\in \WG,$$ with $c_n\in \setC$ be a finite sum of modes
and associated modes with outgoing wavenumbers $\kappa_n(\omega)\in \setC$. Then $\bu$ is bounded and the Laplace
transform in the \dfo{infinite}\dfn{longitudinal} direction \dfo{$\hat \bu(s,\WGy):=$}$\LT(\bu(\bullet,\eta))(s)$
for $(s,\WGy) \in \setR_+\times (-\WGTH,\WGTH)$ takes the form
\begin{align}\label{eq:LTu}
\dfo{\hat\bu(s,\WGy)}\dfn{\LT(\bu(\bullet,\eta))(s)}=\sum_n^N \frac{c_n}{(s-\iu\kappa_n)^{m_n}}\tilde\bw_n(\WGy),
\end{align}
with $m_n\in \setN$. It has a meromorphic extension to $\setC$ with poles at $\{\iu\kappa_n\}$. In contrast, the Laplace 
transform of $e^{-\iu\kappa_n\WGx}$ has a pole at $-\iu\kappa_n$ (see Fig.~\ref{fig:wn3} for a typical situation).
\dfo{Roughly speaking the}\dfn{The} \emph{pole condition} states, that $\bu$ is called outgoing, if \dfo{$\hat\bu$}\dfn{$\LT\bu$}
has no poles in a suitable region of the complex plane. Here, one might chose a point symmetric, smooth
and asymptotically straight boundary curve $\curve$ separating the outgoing
from the incoming poles, \dfn{simply connected} domains $\curve^-$, $\curve^+$ such that
$\setC=\curve^-\dot\cup\curve \dot\cup \curve^+$ and asks \dfo{$\hat\bu$}\dfn{$\LT\bu$} to be holomorphic in $\curve^-$. Typical choices of $\curve$ can be seen in Fig.~\ref{fig:ComplexWNrs}.
The assumption
\begin{align}\label{eq:ass}
\setWN(\omega)\subset\curve^+,
\end{align}
with $\setWN(\omega):=\{\iu \kappa_1(\omega),\iu \kappa_2(\omega),\dots\}$ being the set of outgoing wavenumbers multiplied with $\iu$ for the frequency $\omega$, is thereby essential.
\dfn{To formulate the above in a more compact way, we introduce the
\emph{Hardy space $H^-(\Gamma)$} as the subspace of $L^2(\Gamma)$, such that for each $f\in H^-(\Gamma)$
there exists an in $\Gamma^-$ holomorphic function $f_\mathrm{vol}$ with $f$ being the non-tangential
limit of $f_\mathrm{vol}$. We refer to \cite[App.~A]{HallaHohageNannenSchoeberl:14} for a detailed definition and
properties of such Hardy spaces. We can formulate the pole condition now in the following way
\begin{align}\label{eq:pc}
\LT\bu\in[H^{-}(\curve)\otimes L^2(\tilde\surf)]^2.
\end{align}
}

\dfo{For the following, we refer to \cite[App.~A]{HallaHohageNannenSchoeberl:14} for the details.
\begin{defi}[Hardy space]
\label{def:HardySpace}
Let $\curve=\gamma(\setR)$ be an oriented curve parameterized by a twice continuously differentiable function $\gamma:\setR\to\setC$ of the form $\gamma(z) = z\sigma(|z|)$, $z\in\setR$ with $\sigma:[0,\infty)\to S^1$, and suppose that
$\sup_{z\in\setR}|\gamma'(z)|<\infty$ and there exists $\sigma_{\infty}\in S^1$ such that $\lim_{|z|\to\infty}|\gamma(z)-\sigma_{\infty}z|=0$. 
Then the Hardy space $H^{\pm}(\curve)$ is the set of functions $f_{\rm bd}\in L^2(\curve)$ such that there exists an analytic function $f_{\rm vol}$ in 
\begin{equation}
 \curve^\pm:=\{\gamma(z)\exp(\pm \iu t): z>0, t\in (0,\pi)\}
\end{equation}
and a sequence of rectifiable Jordan curves $C_1,C_2,\dots$ tending to the boundary in $\curve^\pm$ such that the integrals $\int_{C_n} |f_{\rm vol}(s)|^2|ds|$ are uniformly bounded and $f_{\rm vol}(s)\to f_{\rm bd}(t)$ for $s\in\curve^\pm,s\to t$ non-tangentially for almost all $t\in\curve$. 
\end{defi}
The assumptions on $\curve$ ensure, that the Hardy spaces are well-defined Hilbert spaces with the $L^2(\curve)$ scalar product and contain the dense sets $\{1/(s-\lambda), \lambda\in\curve^\mp\}\subset H^{\pm}(\curve)$. Moreover, $\curve$ is point symmetric (i.e. $-\curve=\curve$),  $-\curve^+=\curve^-$, and $\setC=\curve^-\dot\cup\curve \dot\cup \curve^+$. A typical choice of $\curve$ can be seen in Fig.~\ref{fig:ComplexWNrs}. 
\begin{defi}[pole condition]
\label{def:PoleCondition}
 Let $\curve$ fulfill the assumptions of Def.~\ref{def:HardySpace}. A function $\bu\in[H^1_{\rm loc}(W)]^2$ satisfies the pole condition w.r.t. the curve $\curve$ if
\begin{align*}
\int_0^{\infty}e^{-\WGx_0 \WGx}\|\bu(\WGx,\cdot)\|_{[L^2(\tilde\surf)]^2}d\WGx < \infty
\end{align*}
for some $\WGx_0>0$ and the Laplace transform $(\LT\bu)(s):=\int_0^{\infty}e^{-s \WGx}\bu(\WGx,\cdot)d\WGx$ (with values in $[L^2(\tilde\surf)]^2$) has a holomorphic extension from $\{s\in\setC\colon\Re s>\WGx_0\}$ to $\curve^-$ with $[L^2(\tilde\surf)]^2$-boundary values on $\curve$ such that
\begin{align*}
\LT\bu\in[H^{-}(\curve)\otimes L^2(\tilde\surf)]^2.
\end{align*}
\end{defi}
In order to connect the pole condition to the modal radiation condition Def.~\ref{defi:modal_rad_cond}, we need additional
assumptions on $\curve$.
\begin{ass}\label{ass:Gamma}
Let $0\leq\omega_l<\omega_u$ be given. Let $\setWN(\omega):=\{\iu \kappa_1(\omega),\iu \kappa_2(\omega),\dots\}$ be the set of outgoing wavenumbers multiplied with $\iu$ for the frequency $\omega\in\setR^+\setminus\{\omega_1,\omega_2,\dots\}$. We assume, that $\curve$ satisfies the assumptions in Def.~\ref{def:HardySpace} and is chosen such that
\begin{equation}\label{eq:AssumptWaveNumbersGamma}
\bigcup_{\omega\in(\omega_l,\omega_u)\setminus\{\omega_1,\omega_2,\dots\}} \setWN(\omega)\subset \curve^+.
\end{equation}
\end{ass}}
If \dfn{\eqref{eq:ass}}\dfo{Ass.~\ref{ass:Gamma}} holds true, \dfn{then}\dfo{than for the interval of frequencies $(\omega_l,\omega_u)\setminus\{\omega_1,\omega_2,\dots\}$} the modal radiation condition \dfo{Def.~\ref{defi:modal_rad_cond}}\dfn{of Sec.~\ref{Sec:ModalAnalysis}} is equivalent to
the pole condition \dfn{\eqref{eq:pc}}\dfo{Def.~\ref{def:PoleCondition}} for solutions $\bu$ to Equ.~\eqref{eq:modal} having the form
$\bu=\sum_{n=1}^N c_n\bu_n$ and $\bu_n$ being modes and associated modes. 

\dfn{
The pole condition \eqref{eq:pc} does not depend on the wavenumbers but only on $\curve$ and is therefore frequency independent. The wavenumbers are hidden in the assumption \eqref{eq:ass}. Of course, we have to discuss how to chose $\curve$ (see Sec.~\ref{Sec:CondCurve}). Moreover, for the interpretation of solutions to the resonance problem \eqref{eq:res_classic} in Sec.~\ref{sec:resprob} we need the wavenumbers. But the numerical method based on \eqref{eq:pc} presented in the following sections is independent of the wavenumbers and waveguide modes. This facilitates the implementation of the method a lot. Moreover, a discretization of the resonance problem \eqref{eq:res_classic} leads to a linear generalized matrix eigenvalue problem. 
}

\section{\dfn{The Hardy space infinite element}\dfo{Hardy space infinite elements}}\label{sec:HSM}
\dfn{The pole condition $\LT\bu\in[H^{-}(\curve)\otimes L^2(\tilde\surf)]^2$ is a nice way to reformulate
the physical radiation condition, but a stable numerical method based on this framework is delicate.
In \cite{HallaHohageNannenSchoeberl:14} such a method was developed for a one dimensional toy problem including
a complete convergence analysis. We present here only the results needed for an implementation
in the setting of
a convected
one dimensional Helmholtz equation
\begin{subequations}\label{eq:convhh}
\begin{align}
\label{eq:convhh-a} -u''+u'-\omega^2u &=0, \quad x>0,\\
\label{eq:convhh-b} u'(0)&=u_0',\\
\label{eq:convhh-c} \LT u&\in H^-(\curve),
\end{align}
\end{subequations}
whereas the curve $\curve$ is chosen such that the solution to \eqref{eq:convhh} is unique.
Since solutions of \eqref{eq:convhh-a} have the form $C_1e^{\iu\kappa_1}+C_1e^{\iu\kappa_1}$
with $C_1,C_2,\kappa_1,\kappa_2\in\setC$, the curve $\Gamma$ is such that $\iu\kappa_1\in\Gamma^+$ and
$\iu\kappa_2\in\Gamma^-$ (or vice-versa).}

\dfn{The discretization of the elasticity problem discussed in Section~\ref{sec:disc-method}
will be straightforward. In order to enhance readability, we have refrained from giving the correct mathematical
framework including the variational formulation in the Hardy space. For the elastic waveguide
problem this can be deduced along the lines of \cite{HallaHohageNannenSchoeberl:14} and \cite{HohageNannen:15}.}

\subsection{\dfn{Basis functions and infinite element matrices}}
\dfn{Similar to a classical infinite element method we start with the variational form
\begin{equation}\label{eq:1dvarform}
 \int_0^\infty u' v\,dx' + \int_0^\infty u'v \,dx -\omega^2 \int_0^\infty u v \, dx = -u_0' v(0)
\end{equation}
of \eqref{eq:convhh}, which holds true for the solution $u$ and sufficiently fast decaying, smooth test functions $v$.
We are using a Galerkin scheme with the same basis functions $\bax_j$, $j\in \setN$, as ansatz and test functions. These basis functions should fulfill the
pole condition $\LT \bax_j \in H^-(\curve)$.}

\dfn{From \cite{HallaHohageNannenSchoeberl:14} we know that for  any two complex parameters $\pol,\poli \in \curve_+$
the linear hull of
\begin{equation}\label{eq:basisHardyspace}
  \psi_j^{\pol,\poli}(s):=\frac{\pol+\poli}{s-\poli}\left(\frac{s+\pol}{s-\pol}\right)^{\lfloor (j+1)/2 \rfloor}
\left(\frac{s+\poli}{s-\poli}\right)^{\lfloor j/2 \rfloor}, \quad j\in\setN_0,
\end{equation}
is dense in the Hardy space $H^-(\curve)$. The basis functions $\bax_j$ of our Galerkin scheme are defined via their Laplace transforms 
\begin{equation}
 (\LT \bax_1)(s):=\frac{1}{s-\pol},\qquad (\LT \bax_j)(s):= \frac{\psi_{j-2}^{\pol,\poli}(s)}{s-\pol},\quad j=2,\dots,\qquad s \in \curve.
\end{equation}
They have the form
$$\bax_j(x)=e^{\pol x}p_{j}^{\pol}(x) + e^{\poli x}p_{j}^{\poli}(x),\qquad x \geq 0,$$
with polynomials $p_{j}^{\pol},p_{j}^{\poli}$ and in particular $\spann\{\bax_0,\bax_1\}=\spann\{e^{\pol x},e^{\poli x}\}$. Comparing the basis functions with the modes \eqref{eq:mode}
shows, that these functions represent the longitudinal part of the modes exactly for wavenumbers $\kappa=\pol/\iu$ and $\kappa=\poli/\iu$. Therefore, $\pol/\iu$ and $\poli/\iu$ can be considered as two different wavenumbers. The actual choice of these parameters will be discussed in the next subsections.}

\dfn{It is easy to show that $\bax_j$, $j\in \setN$ fulfill the pole condition and moreover, by a limit theorem of the Laplace transform there holds
$$ \bax_j (0)=\lim_{x \searrow 0} \bax_j (x)=\lim_{s \to \infty} s (\LT \bax_j)(s)
 =\begin{cases}1,&\quad j=1\\ 0,&\quad j>1\end{cases}.$$
This facilitates coupling of these basis functions at $x=0$ with e.g. standard finite element basis functions for $x<0$.} 

\dfn{We are left with the computation of the integrals in \eqref{eq:1dvarform} using $ \bax_j$ for $u$ and $v$. If $\pol$ and $\poli$ are chosen with negative real part,
the integrals are bounded. Similar to \cite[Lemma A.1]{HohageNannen:09} it can be shown, that
\begin{equation}
\int_0^\infty \bax_j(x) \bax_k(x) \, dx = \frac{-\iu}{2 \pi} \int_\curve (\LT \bax_j)(s) (\LT \bax_k)(-s),\,ds,\qquad j,k \in \setN.
\end{equation}
 Since $\LT \bax_j$ are meromorphic functions, the integrals on the right hand side can be computed by the residue theorem. Following the computations in \cite{HallaHohageNannenSchoeberl:14}, the infinite mass matrix 
 $\left(\Mx^\infty\right)_{j,k}:=\int_0^\infty  \bax_j(x) \bax_k(x) \, dx$ is given by the tridiagonal matrix
 \begin{equation}
  \Mx^\infty = \frac{-1}{\pol\poli} \smat{ \boxed{\begin{smallmatrix} \poli & -\poli \\ -\poli & \pol+\poli\end{smallmatrix}}  & \boxed{\begin{smallmatrix}  \hspace{7pt} 0  \hspace{7pt} &  \hspace{7pt} 0  \hspace{7pt} \\ -\pol & 0\end{smallmatrix}} & &\\
   \boxed{\begin{smallmatrix}  \hspace{4pt} 0  \hspace{4pt} &  \hspace{4pt} -\pol  \hspace{4pt} \\ 0 & 0\end{smallmatrix}}  & \boxed{\begin{smallmatrix} \pol+\poli & -\poli \\ -\poli & \pol+\poli\end{smallmatrix}} & \ddots  & \\
    & \ddots & \ddots}.
 \end{equation}
Analogue calculations can be performed for $\partial_x \bax_j$, since
 \begin{equation*}
  (\LT \partial_x \bax_1)(s)=\frac{\pol}{s-\pol},\qquad  (\LT \partial_x \bax_j)(s)= \frac{s \, \psi_{j-2}^{\pol,\poli}(s)}{s-\pol},\quad j=2,\dots.
 \end{equation*}
The drift $\left(\Dx^\infty\right)_{jk}:=\int_0^\infty  \partial_x \bax_j(x) \bax_k(x) \, dx$ and stiffness $\left(\Sx^\infty\right)_{jk}:=\int_0^\infty  \partial_x \bax_j(x) \partial_x \bax_k(x) \, dx$
matrices are given by
\begin{equation}
\Dx^\infty = \frac{1}{2} \smat{\phantom{-}1&\phantom{-}1&\phantom{\ddots}&\phantom{\ddots}&\phantom{\ddots}\\-1&\phantom{-}0&\phantom{-}1&\phantom{\ddots}\\\phantom{\ddots}&-1&\phantom{-}0&\phantom{-}\ddots\\\phantom{\ddots}&\phantom{\ddots}&\phantom{-}\ddots&\phantom{-}\ddots
 },\qquad 
 \Sx^\infty = \frac{-1}{2} \smat{ \boxed{\begin{smallmatrix} \pol & \pol \\ \pol & \pol+\poli\end{smallmatrix}}  & \boxed{\begin{smallmatrix}  \hspace{7pt} 0  \hspace{7pt} &  \hspace{7pt} 0  \hspace{7pt} \\ \poli & 0\end{smallmatrix}} & &\\
   \boxed{\begin{smallmatrix}  \hspace{4pt} 0  \hspace{4pt} &  \hspace{4pt} \poli  \hspace{4pt} \\ 0 & 0\end{smallmatrix}}  & \boxed{\begin{smallmatrix} \pol+\poli & \pol \\ \pol & \pol+\poli\end{smallmatrix}} & \ddots  & \\
    & \ddots & \ddots}.
\end{equation}}
\dfn{Using only the first $N^{\rm long}$ basis functions and therefore only the $N^{\rm long} \times N^{\rm long}$ upper left block of these matrices leads to the discretization of \eqref{eq:1dvarform}: Find $U^{(N^{\rm long})}\in \setC^{N^{\rm long}}$ such that
\begin{equation}\label{eq:1dGLS}
 \left(\Sx^{N^{\rm long}} + \Dx^{N^{\rm long}} -\omega^2 \Mx^{N^{\rm long}}\right) U^{(N^{\rm long})} = (-u_0',0,\dots,0)^\top.
\end{equation}
Note, that the matrices $\Sx^{N^{\rm long}}$, $\Dx^{N^{\rm long}}$ and $\Mx^{N^{\rm long}}$ are independent of the angular frequency $\omega$.}

\subsection{\dfn{Relation between the method parameters $\pol$, $\poli$ and $\curve$}}\label{Sec:defCurve}
\dfn{\eqref{eq:1dGLS} is a conforming discretization of \eqref{eq:convhh-a}-\eqref{eq:convhh-b}.
The side constraint $\LT\sum_{j=1}^{N^\textrm{long}}U^{(N^\textrm{long})}_j\bax_j\in H^-(\Gamma)$ (see \eqref{eq:convhh-c})
is however not only fulfilled for the chosen $\Gamma$, but for any 
$\tilde\Gamma$ with $\pol,\poli\in\tilde\Gamma^+$.
If for such a second $\tilde \curve$  the wavenumbers $\kappa_{1,2}$ of \eqref{eq:convhh-a} multiplied with $\iu$ belong to different sides of the curves, e.g.
$\iu\kappa_1\in\Gamma^+\cap\tilde\Gamma^-$ and $\iu\kappa_2\in\Gamma^-\cap\tilde\Gamma^+$,
Equations \eqref{eq:convhh-a}-\eqref{eq:convhh-b} together with the condition $\LT u\in H^-(\tilde\Gamma)$ define a different solution then \eqref{eq:convhh}. But \eqref{eq:1dGLS} is again a conforming discretization for this second problem. This leads to the question that if the solutions
$\sum_{j=1}^{N^{\rm long}} U^{(N^{\rm long})}_j \bax_j$ of \eqref{eq:1dGLS} converge (in any norm)
to a function $U$, does $\LT U\in H^-(\Gamma)$ or $\LT U\in H^-(\tilde\Gamma)$ hold?}

\dfn{The answer was given in \cite{HallaHohageNannenSchoeberl:14}: The linear hull of the functions
$\LT\bax_j,j\in\setN$, is dense in any $H^-(\Gamma)$ with $\pol,\poli\in\Gamma^+$.
But they form a stable basis only in one specific space $H^-(\Gamma_{\pol,\pol})$,
i.e.~any $U\in H^-(\Gamma_{\pol,\poli})$ can be expanded $U=\sum_{j\in\setN}\alpha_j\LT\bax_j$
with a square summable sequence $(\alpha_j)_{j\in\setN}$.
One can deduce that if $\{\iu\kappa_1,\iu\kappa_2\}\cap\Gamma^\pm=
\{\iu\kappa_1,\iu\kappa_2\}\cap\Gamma_{\pol,\poli}^\pm$
Eq.~\eqref{eq:1dGLS} is uniquely solvable at least for sufficiently large $N^\textrm{long}$
with solutions $U^{(N^{\rm long})}$ and
$\lim_{N^\textrm{long}\to\infty}\|\LT u-\sum_{j=1}^{N^{\rm long}} U^{(N^{\rm long})}_j\LT\bax_j\|_{L^2(\Gamma_{\pol,\poli})}=0$ holds.
Thence the choice of $\pol,\poli$ for \eqref{eq:1dGLS} implicitly poses
the pole condition w.r.t.~$\Gamma_{\pol,\poli}$ on the solution $u$ of \eqref{eq:convhh}.}

\dfn{The curve $\curve_{\pol,\poli}$ is given by the algebraic variety
\begin{equation}
\label{eq:defcurve}
\curve_{\pol,\poli}:=\left\{s\in\setC\colon\left|\frac{s+\pol}{s-\pol}\,\frac{s+\poli}{s-\poli}\right|=1\right\}.
\end{equation}
In order to have a physically correct pole condition for elastic waveguide problems 
we are left with the task to find $\pol,\poli$ such that $\mathbb{S}(\omega)\subset\Gamma_{\pol,\poli}^+$,
which will be the issue of the next subsection.
}
\dfn{The following properties of the curve $\curve_{\pol,\poli}$ can all be deduced from \eqref{eq:defcurve}: 
First of all for $g_{\pol,\poli}(s):= \frac{|s-\pol|}{|s+\pol|}\frac{|s-\poli|}{|s+\poli|}$ we have the characterizations 
\begin{equation}
 \Gamma^+_{\pol,\poli}=\{s\in\setC\colon  g_{\pol,\poli}(s)<1\},\quad \Gamma_{\pol,\poli}=\{s\in\setC\colon  g_{\pol,\poli}(s)=1\},\quad \text{and } \Gamma^-_{\pol,\poli}=\{s\in\setC\colon  g_{\pol,\poli}(s)>1\}.
\end{equation}
The curve can explicitly be parameterized by $\curve_{\pol,\poli}=\gamma_{\pol,\poli}(\setR)$ with 
\begin{equation}\label{eq:ParamCurve}
 \gamma_{\pol,\poli}(r):=-\mathfrak{i}r\frac{r^2(\pol+\poli)+|\pol|^2\poli+|\poli|^2\pol}{|r^2(\pol+\poli)+|\pol|^2\poli+|\poli|^2\pol|},\qquad r\in\setR,
\end{equation}
which shows that its asymptotic behavior is a line. In the case $|\pol|=|\poli|$ it actually is a line.
Further scaling the poles scales the curve: $\Gamma_{r\pol,r\poli}=r\Gamma_{\pol,\poli}$ for all $r>0$.
If we chose $\pol,\poli$ such that
\begin{subequations}\label{eq:conds0s1}
 \begin{align}
  \label{eq:poleq1} \Re(\pol), \Re(\poli) &<0,\\
  \label{eq:poleq2} \Im(\pol+\poli) &>0,\\
  \label{eq:poleq3} |\pol|^2\Im(\poli)+|\poli|^2\Im(\pol) &<0,
 \end{align}
\end{subequations}
it holds $\iu(-\mpara(s_0,s_1),0)\cup\iu(\mpara(s_0,s_1),\infty)\subset\Gamma_{\pol,\poli}^+$ and
$\iu(0,\mpara(s_0,s_1))\cup\iu(-\infty,-\mpara(s_0,s_1))\subset\Gamma_{\pol,\poli}^-$
with
\begin{align}
\mpara(s_0,s_1):=\sqrt{-\tfrac{|\pol|^2\Im(\poli)+|\poli|^2\Im(\pol)}{\Im(\pol+\poli)}}.
\end{align}}

\subsection{Choice of $\pol$ and $\poli$}\label{Sec:CondCurve}
\dfn{Coming back to the elastic waveguide problem, we inspect the location and qualitative properties of the wavenumber spectrum in order to find reasonable values for
the parameters $\pol$ and $\poli$.}
First, let us look only on the real wavenumbers in Fig.~\ref{fig:disprelfrequencyrange} \dfn{for fixed parameters $\nu$
and $H,\rho,E$}.
There exist eight frequencies $\omega>0$ \dfn{which support waves} with vanishing group velocities
$\left(\frac{\partial \kappa_n}{\partial \omega}\right)^{-1}$. 
Six of them yield vanishing wavenumbers $\kappa_n(\omega)=0$ and two don't.
These $8$ frequencies $\omega_1,\dots,\omega_8$ are part of the sequence of frequencies defined in Rem.~\ref{Rem:TheoModes}
and cannot be treated by the pole condition, since for $\omega \to \omega_j$
for $j=1,\dots,8$ an outgoing wavenumber converges to an incoming wavenumber.
Hence, outgoing wavenumbers cannot be separated by $\curve$ from incoming wavenumbers.
In Fig.~\ref{fig:disprelfrequencyrange} we have made a (non unique) decomposition into five intervals $(a,b)$ such that 
for all $\omega\in (a,b)\setminus\{\omega_1,\dots,\omega_8\}$ one of the two following cases hold true:
\begin{figure}
 \centering
\includegraphics[width=\textwidth]{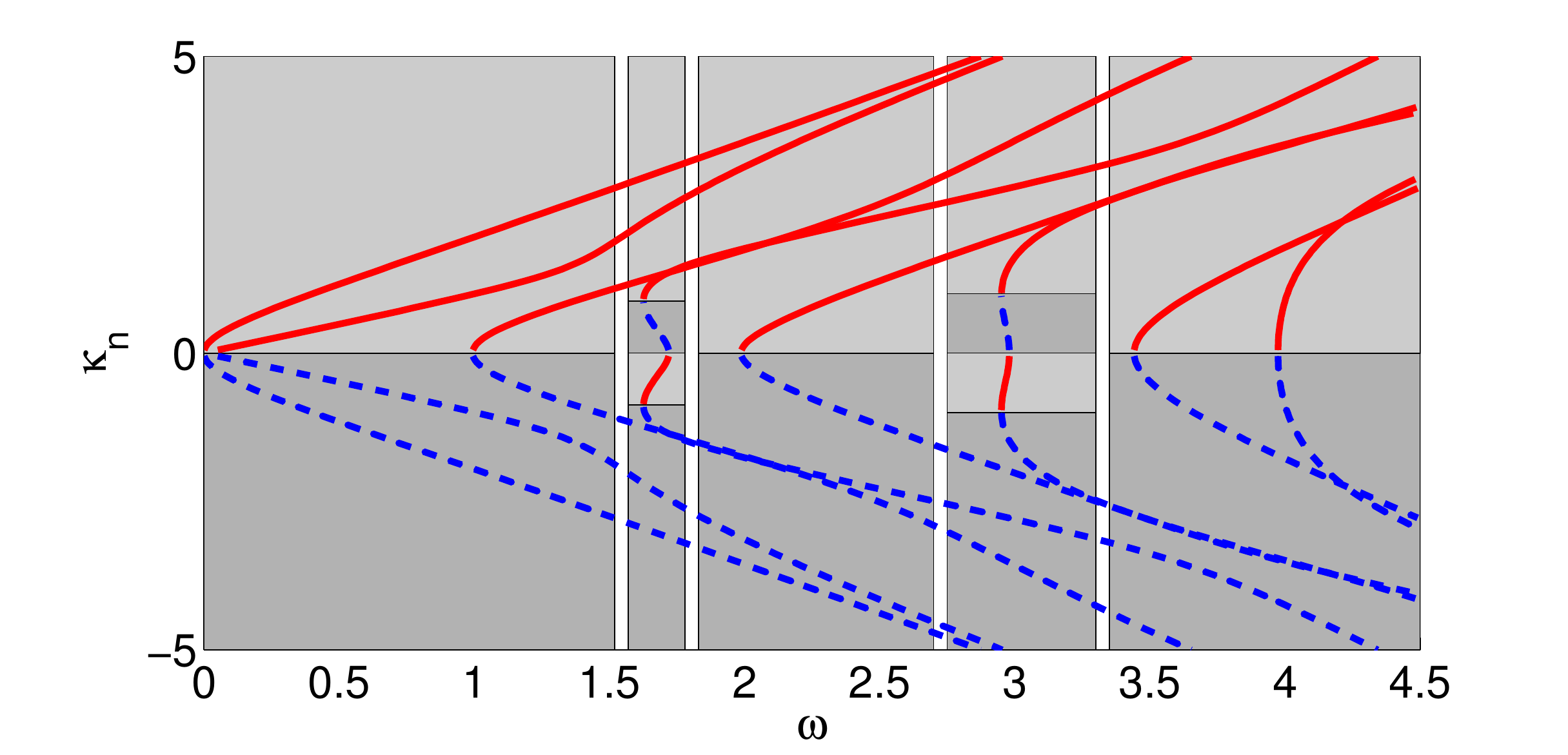}
 \caption{The first eight dispersion curves for $H=\rho=E=1$ and $\nu=0.25$.  
Modes corresponding to the red solid part have positive group velocity and modes corresponding to the blue dashed
part have negative group velocity.}
 \label{fig:disprelfrequencyrange}
\end{figure}
\begin{enumerate}
 \item For all real wavenumbers the signs of group velocity and phase velocity coincide. Hence, the real wavenumbers 
 with positive group velocity belong to $(0,\infty)$.
 \item There exists $\theta>0$, such that the real wavenumbers with positive group velocity belong to $(-\theta,0)\cup (\theta,\infty)$.
\end{enumerate}
Hence, the first case requires $\iu(0,\infty) \subset\curve^+$ and
the second case $\iu(-\theta,0)\cup \iu(\theta,\infty)\subset\curve^+$.
\dfn{Chosing different parameters $H,\rho,E$ would only result in a different scaling of Fig.~\ref{fig:disprelfrequencyrange}.
On the other hand varying $\nu$ changes the qualitative behavior. Nevertheless, for all
$\omega$ either the first case or the second case seems to be true for cylindrical elastic waveguide problems with constant $E$, $\rho,$ and  $\nu$ in the waveguides
(at least the authors were not able to find any contradicting parameters).}

\dfn{
For the non-real outgoing wavenumbers it is known (see Rem~\ref{Rem:TheoModes}),
that for a fixed frequency they are contained in a circular proper subsector
of $\{s\in\setC\colon\Im s>0\}$.
Because the wavenumbers are continuous with respect to $\omega$,
the above also holds uniformly for $\omega$ in suitable small intervals.}

\begin{figure}
 \centering
\subfigure[\label{fig:ComplexWNrsHP} $\omega \in (0.8,1.57)$, $\pol=\poli=-1+0.2\iu$]{\includegraphics[width=0.49\textwidth]{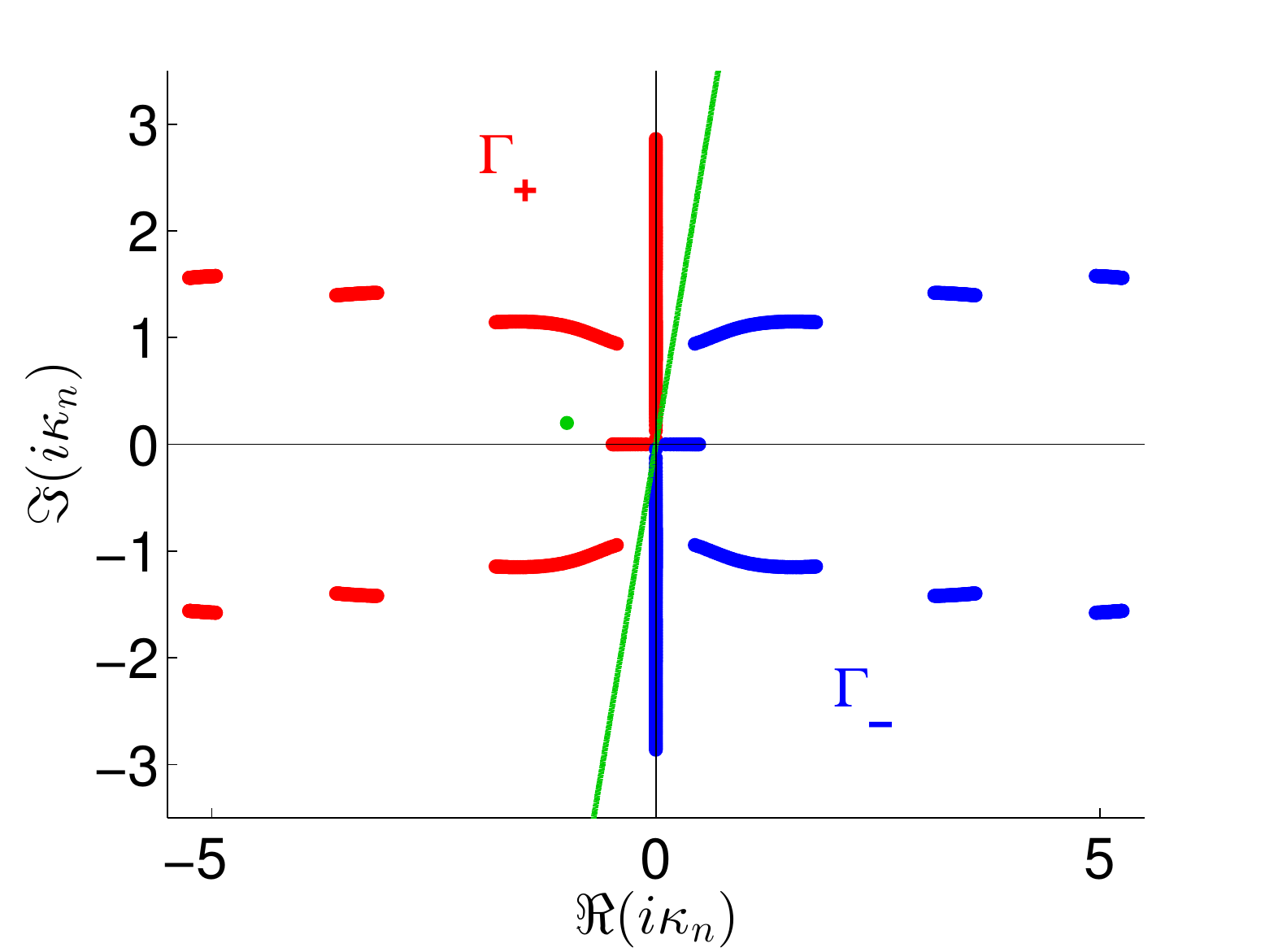}} \hfill
\subfigure[\label{fig:ComplexWNrsCurved}$\omega \in (1.57,1.78)$, $\pol=-0.3742-0.4886\iu$, $\poli=-0.7752+1.0396\iu$]{\includegraphics[width=0.49\textwidth]{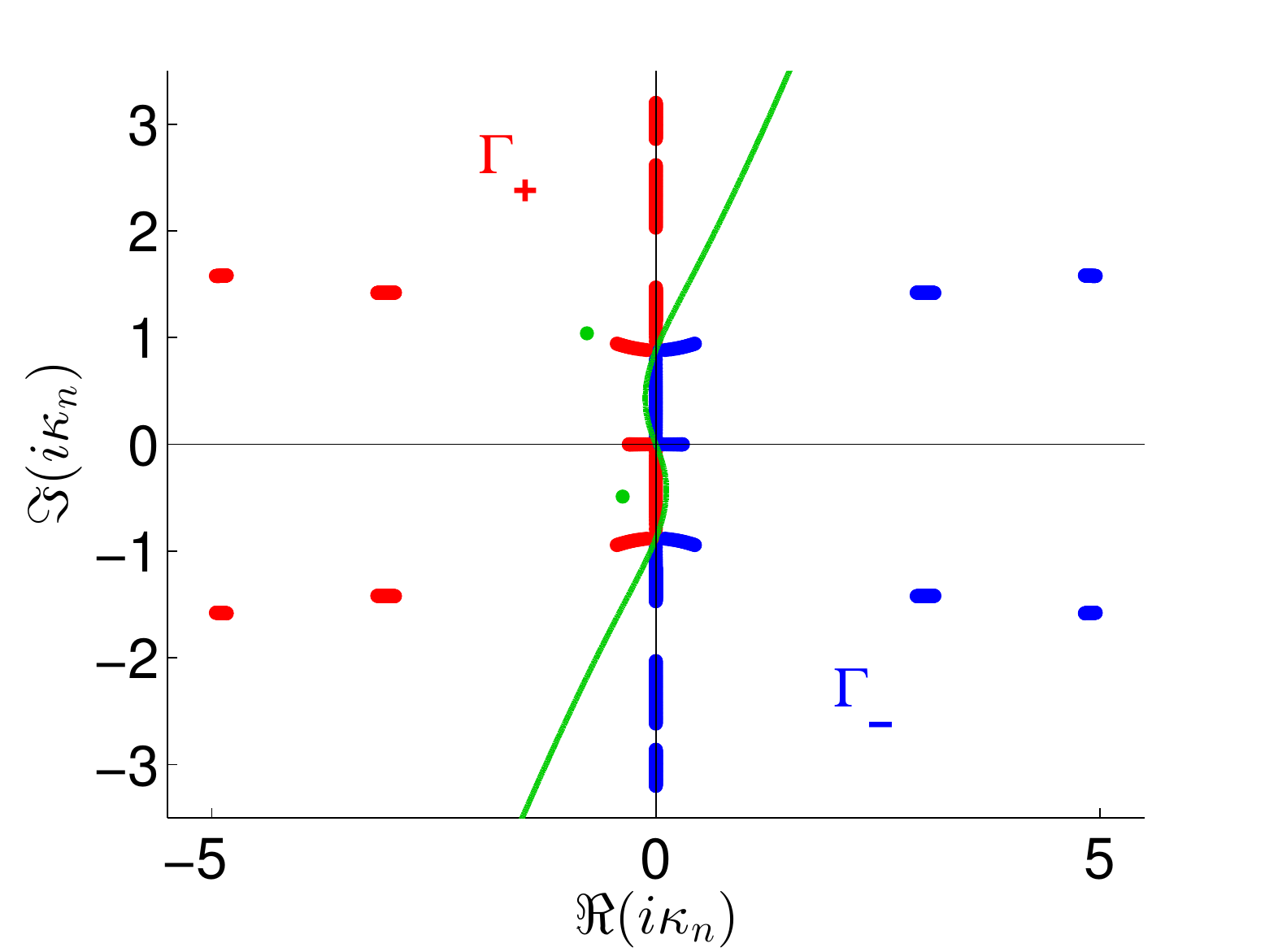}}
 \caption{Lowest outgoing (red) and incoming (blue) wavenumbers (multiplied with $\iu$) for different intervals of 
 frequencies. Separating curves $\curve$ with parameters $\pol$ (and $\poli$) are marked green. Material parameters are $H=\rho=E=1$ and $\nu=0.25$.}
 \label{fig:ComplexWNrs}
\end{figure}

\dfn{For an implementation we have to discuss how to chose the parameters $\pol$ and $\poli$ in order to yield $\mathbb{S}(\omega)\subset\curve^+_{\pol,\poli}$.
For the two cases from above (see Fig.~\ref{fig:ComplexWNrs}) we can proceed as follows:
\begin{enumerate}
 \item $\iu(0,\infty) \subset\curve^+$ (Fig.~\ref{fig:ComplexWNrsHP}):  The choice $|\pol|=|\poli|$ is sufficient. In particular we can chose $\pol=\poli$. Then $\curve_{\pol,\poli}$ is a straight
line $\iu(\pol+\poli)\setR$. Choosing $\Re\big((\pol+\poli)/|\pol+\poli|\big)$ small enough yields
$\mathbb{S}(\omega)\subset\curve^+_{\pol,\poli}$.
\item $\iu(-\theta,0)\cup \iu(\theta,\infty)\subset\curve^+$ (Fig.~\ref{fig:ComplexWNrsCurved}): We chose arbitrary $\tilde\pol,\tilde\poli$, such that \eqref{eq:conds0s1}
is satisfied. E.g.~we can chose $\tilde\pol,\tilde\poli$ as in Fig.~\ref{fig:ComplexWNrsCurved}. Then we scale
$\pol:=\tilde\pol\theta/\mpara(\tilde s_0,\tilde s_1)$
and
$\poli:=\tilde\poli\theta/\mpara(\tilde s_0,\tilde s_1)$,
such that $\iu\setR\cap\curve_{\pol,\poli}=\{0,\pm\iu\theta\}$.
The real valued outgoing wavenumbers are now contained in $-\iu\curve_{\pol,\poli}^+$.
To check if the non-real outgoing wavenumbers are contained in $-\iu\curve_{\pol,\poli}^+$,
$g_{\pol,\poli}(\iu\kappa_n)<1$
can be evaluated. If not, define
$\pol^t:=\pol\theta/\mpara(\pol,t\overline{\pol}+(1-t)\poli)$ and
$\poli^t:=(t\overline{\pol}+(1-t)\poli)\theta/\mpara(\pol,t\overline{\pol}+(1-t)\poli)$
for $t\in[0,1]$. There hold $s_0^0=s_0, s_1^0=s_1$ and $s_0^1=s_0, s_1^1=\overline{s_0}$,
i.e.~$t\mapsto\Gamma_{s_0^t,s_1^t}$ is a homotopy between $\Gamma_{s_0,s_1}$ and $\iu\setR$,
such that $\Gamma_{s_0^t,s_1^t}\cap\iu\setR=\{0,\pm\iu\theta\}$ for all $t\in[0,1)$.
Choosing $t\in[0,1]$ close enough to one, $\mathbb{S}(\omega)\subset\curve^+_{\pol^t,\poli^t}$
can be ensured.
\end{enumerate}
The approximation error of the method depends on $g_{\pol,\poli}(\iu \kappa_n)^{N_{\rm long}}$. Hence, in general one should try to find $\pol$ and $\poli$ such that $g_{\pol,\poli}(\iu \kappa_n)$ 
is as small as possible for the first wavenumbers $\kappa_n$.
}

\dfn{Note, that the wavenumbers and modes are not used in the numerical method. A rough estimate of the location of the wavenumbers is needed in order to chose reasonable values for $\pol$ and $\poli$, but the 
Hardy space infinite element method itself is independent of the wavenumbers and the waveguide modes.}

\section{\dfn{Tensor product discretization of one waveguide}\dfo{Discretization method}}\label{sec:disc-method}
\dfo{We do not give a convergence analysis in this paper, but merely report our numerical method and experiments. Thus we omit to build a variational framework and present our discretization of~\eqref{eq:lame} in a straightforward way. The authors hope that this enhances the readability.}
\dfn{For a discretization of \eqref{eq:lame} a standard finite element method for the bounded interior domain $\Oi$ can be used. For the waveguides a tensor product method appears reasonable. We present here only shortly such a tensor product method for one reference waveguide. The extension to several waveguides and the coupling with the finite element method for $\Oi$ is straightforward.}

\dfo{
Multiplying \eqref{eq:lame-inner} with test functions $\bv$ and partial integration over $\Oi$ leads to the variational formulation for the interior problem
\begin{subequations}
\begin{equation}\label{eq:varformint}
  \Bai{\bui}{\bvi} - \omega^2 \Bbi{\bui}{\bvi} =\Lf{\bvi} + \Baux{\bui}{\bvi}
\end{equation}
with
\begin{eqnarray}
\Bai{\bui}{\bvi}&:=&\int_{\Oi} \left(2 \mu \bepsilon (\bui):\bepsilon (\bvi) + \lambda \dv \bui \, \dv \bvi\right) d\bx,\\
\Bbi{\bui}{\bvi}&:=&\int_{\Oi} \rho \, \bui \cdot \bvi\,d\bx,\\
 \Lf{\bvi}&:=&\int_{\Oi} \bfc \cdot \bvi \, d\bx + \int_{\partial \Oi \setminus \left(\bigcup_{\ell=1}^L \surf_\ell \right)} \bg \cdot \bvi \, ds,\\
 \Baux{\bui}{\bvi}&:=&\sum_{\ell=1}^\nrWG \int_{\surf_\ell} \left(\bsigma(\bui) \bn\right) \cdot \bvi \, ds.
\end{eqnarray} 
\end{subequations}
Here, we have used Neumann boundary conditions at $\partial \Oi \setminus \left(\bigcup_{\ell=1}^L \surf_\ell \right)$ and the notation $A:B=\sum_{i,j=1}^2 A_{ij} B_{ij}$ for $A,B \in \setC^{2 \times 2}$. The term $\Baux{\bui}{\bvi}$ will cancel out later on with the boundary terms of the variational formulations in the waveguides.}

\dfo{We are using for the bounded domain $\Oi$ continuous, high order triangular finite elements in $H^1(\Oi)$ based on integrated Jacobi polynomials (see e.g. \cite{SchoeberlZaglmayr:05,Zaglmayr:06}). 
If we denote the scalar basis functions by $\bai_j$ for $j=1,\dots N^{\rm int}$ the left hand side of \eqref{eq:varformint} leads to the system matrix 
$\Mai - \omega^2 \Mbi \in \setC^{2 N^{\rm int} \times 2 N^{\rm int}}$ with 
\begin{equation}
\begin{aligned}
 \Mai&:=\smat {\Bai{\smat{\bai_j\\ 0}}{\smat{\bai_k \\ 0}}_{j,k=1}^{N^{\rm int}} &  \Bai{\smat{\bai_j\\ 0}}{\smat{0 \\ \bai_k }}_{j,k=1}^{N^{\rm int}} \\
 \Bai{\smat{0\\\bai_j}}{\smat{\bai_k \\ 0}}_{j,k=1}^{N^{\rm int}} &  \Bai{\smat{0\\ \bai_j}}{\smat{0 \\ \bai_k }}_{j,k=1}^{N^{\rm int}}}, \qquad 
  \Mbi:=\smat {\Bbi{\smat{\bai_j\\ 0}}{\smat{\bai_k \\ 0}}_{j,k=1}^{N^{\rm int}} &  \Bbi{\smat{\bai_j\\ 0}}{\smat{0 \\ \bai_k }}_{j,k=1}^{N^{\rm int}} \\
 \Bbi{\smat{0\\\bai_j}}{\smat{\bai_k \\ 0}}_{j,k=1}^{N^{\rm int}} &  \Bbi{\smat{0\\ \bai_j}}{\smat{0 \\ \bai_k }}_{j,k=1}^{N^{\rm int}}}.
\end{aligned}
 \end{equation}
The integration is done by Gauss quadrature.}

\dfo{Let us assume for this subsection,  that there exists only one waveguide of the form $\WG:=\setR^+\times \tilde \surf$, i.e. $\Omega=\Oi \cup \WG \cup \surf$.} 
\dfn{We consider again the reference waveguide $\WG:=\setR^+\times (-\WGTH,\WGTH)$}.  
For sufficiently fast decaying test functions $\bv$ we obtain as in the last section a variational formulation of \eqref{eq:modal} with the pole condition as radiation condition: \dfo{The bilinear
forms are the same as in the last section, but due to the assumption, that $\bfc$ and $\bg$ have compact support in $\partial\Oi\cap\partial\Omega$ and  $\Oi$, there is no linear form. In full detail, we have}
\dfn{We are looking for $\bu$ with $\LT \bu\in [H^{-}(\curve_{\pol,\poli})\otimes L^2(-\WGTH,\WGTH)]^2$ such that for all suitable $\bv$}
\begin{subequations}
\begin{equation}\label{eq:varformext}
 \Bae{\bu}{\bv} - \omega^2 \Bbe{\bu}{\bv} = - \int_{\surf} \left(\bsigma(\bu) \smat{1\\0} \right) \cdot \bv \, ds
\end{equation}
\dfn{with $\Bae{\bue}{\bve}:=\int_{\WG} \left(2 \mu \,\bepsilon (\bu):\bepsilon (\bv) + \lambda \dv \bu \, \dv \bv\right)  \,d(\WGx,\WGy)$ and $\Bbe{\bu}{\bv}:=\int_{\WG} \rho \,\bu \cdot \bv \, d(\WGx,\WGy)$. In full detail we have}
\begin{eqnarray}\label{eq:BilFaWG}
&&\begin{aligned}
 \Bae{\bu}{\bv}=& \int_{\WG} \left(\left(2\mu+\lambda \right) \partial_\WGx \bu_1 \, \partial_\WGx \bv_1 + \mu \partial_\WGy \bu_1 \, \partial_\WGy \bv_1 \right) \,d(\WGx,\WGy)+\int_{\WG} \left(\left(2\mu+\lambda \right) \partial_\WGy \bu_2 \, \partial_\WGy \bv_2 + \mu \partial_\WGx \bu_2 \, \partial_\WGx \bv_2 \right) \,d(\WGx,\WGy)\\
 &+\int_{\WG} \left(\mu \partial_\WGy \bu_1\, \partial_\WGx \bv_2 + \lambda \partial_\WGx \bu_1\, \partial_\WGy \bv_2 \right) \,d(\WGx,\WGy)+\int_{\WG} \left(\mu \partial_\WGx \bu_2\, \partial_\WGy \bv_1 + \lambda \partial_\WGy \bu_2\, \partial_\WGx \bv_1 \right) \,d(\WGx,\WGy),
\end{aligned}\\
&&\Bbe{\bu}{\bv}:=\int_{\WG} \rho \left(\bu_1\,\bv_1+ \bu_2\,\bv_2  \right) \,d(\WGx,\WGy).
\end{eqnarray}
\end{subequations}
\dfo{Since for a solution $\bu$ to \eqref{eq:lame} and $\bui:=\bu|_{\Oi}$, $\bue:=\bu|_\WG$  we have $\bsigma(\bue) \smat{1\\0} = \bsigma(\bui) \smat{1\\0}$, the right hand side of \eqref{eq:varformext} cancels out with the corresponding term in \dfo{$\Baux{\bui}{\bvi}$} if the test functions are chosen such that $\bve|_\surf=\bvi|_\surf$.}
\dfn{The right hand side of \eqref{eq:varformext} cancels out with the corresponding term of the variational formulation of the interior problem, if 
the test functions are chosen to be continuous on the interface $\surf:=\{0\}\times (-\WGTH,\WGTH)$.}

\dfo{We consider the waveguide} \dfn{The whole waveguide is considered}  as an (in)finite element with element matrix $\Mae -\omega^2 \Mbe$ and tensor product basis functions of the form $\bax_j\otimes \bay_l$. More precisely, for $\bu$ we use the ansatz
\begin{equation}\label{eq:Approxu}
 \bu (\WGx,\WGy)\approx \bu_{N^{\rm long},N^{\rm trans}} (\WGx,\WGy):=\sum_{j=1}^{N^{\rm long}} \sum_{l=1}^{N^{\rm trans}} \bax_j(\WGx)\, \bay_l(\WGy) \smat{\alpha_{jl}^{(1)} \\ \alpha_{jl}^{(2)}},\, \qquad
 (\WGx,\WGy) \in (0,\infty)\times(-\WGTH,\WGTH)
\end{equation}
with $\alpha_{jl}^{(1)},\alpha_{jl}^{(2)}\in \setC$. 
 
In order to ensure continuity of the discrete solution at the interface $\surf$ the basis functions $\bay_l$ are the non-vanishing traces of the interior basis functions  $\bai_j \in H^1(\Oi)$, i.e. for $\bai_{j(1)}|_{\surf},\dots,\bai_{j({N^{\rm trans}})}|_{\surf}\neq 0$
\begin{equation*}
 \bay_l(\WGy):=\bai_{j(l)} (0,\WGy),\qquad \WGy \in (-\WGTH,\WGTH),\qquad l=1,\dots,N^{\rm trans}.
\end{equation*}
The surface matrices $\My,\Dy,\Sy \in \setC^{N^{\rm trans}\times N^{\rm trans}}$ defined by
\begin{equation}
\begin{aligned}
\My :=\left(\int_{-\WGTH}^{\WGTH} \bay_l \, \bay_m\, d\WGy\right)_{l,m=1}^{N^{\rm trans}},\quad \Dy :=\left(\int_{-\WGTH}^{\WGTH} \partial_\WGy \bay_l \, \bay_m\, d\WGy\right)_{l,m=1}^{N^{\rm trans}},\quad
\Sy :=\left(\int_{-\WGTH}^{\WGTH} \partial_\WGy  \bay_l \, \partial_\WGy \bay_m\, d\WGy\right)_{l,m=1}^{N^{\rm trans}}
\end{aligned}
\end{equation}
can be computed e.g. by Gauss quadrature. 

\dfo{We will give later on the basis functions $\bax_1,\dots,\bax_{N^{\rm long}}$ in the infinite direction. Let us assume for the moment that $\bax_j$ fulfill some basic properties, i.e. the pole condition $\LT\bax_j\in H^-(\curve_{\pol,\poli})$, a separation in a boundary basis function $\bax_1$ and inner basis functions $\bax_2,\dots,\bax_{N^{\rm long}}$, i.e.  
\begin{equation}
 \label{eq:BdassumptBasisx}
 \bax_1(0)=1, \qquad \bax_j(0)=0,\quad j=2,\dots,N^{\rm long},
\end{equation}
and that we have the element matrices
\begin{equation}\label{eq:RadMatrices}
\begin{aligned}
\Mx :=\left(\int_0^\infty \bax_j \, \bax_k\, d\WGx\right)_{j,k=1}^{N^{\rm long}},\quad \Dx :=\left(\int_0^\infty \partial_\WGx \bax_j \, \bax_k\, d\WGx\right)_{j,k=1}^{N^{\rm long}},\quad
\Sx :=\left(\int_0^\infty \partial_\WGx \bax_j \, \partial_\WGx  \bax_k\, d\WGx\right)_{j,k=1}^{N^{\rm long}}
\end{aligned}
\end{equation}
in infinite direction at hand.}
\dfn{For the unbounded $\WGx$-direction we use the basis functions and element matrices $\Mx$, $\Dx$ and $\Sx$ of the last section (we have droped the superscripts for convenience). }
The system matrix of \eqref{eq:varformext} is therefore given by $\Mae -\omega^2 \Mbe \in \setC^{2 N^{\rm long}N^{\rm trans}\times 2 N^{\rm long}N^{\rm trans}}$ with
\begin{subequations}
\label{eq:IEMatrices}
\begin{align}
 \Mae&:=\begin{pmatrix}
  (2\mu+\lambda) \Sx \otimes \My + \mu \Mx \otimes \Sy & \mu\, \Dx^\top \otimes \Dy + \lambda \, \Dx \otimes \Dy^\top\\ 
  \mu\, \Dx \otimes \Dy^\top + \lambda \, \Dx^\top \otimes \Dy & (2\mu+\lambda) \Mx \otimes \Sy + \mu \Sx \otimes \My  
 \end{pmatrix},\\
 \Mbe&:=\begin{pmatrix}
         \rho \, \Mx \otimes \My & \bf{0}\\
         \bf{0} & \rho \, \Mx \otimes \My
         \end{pmatrix}.
\end{align}
\end{subequations}
\dfo{Note, that we are working with bilinear forms. Therefore, $A^\top$ with $A\in \setC^{n\times n}$ denotes the transpose
of $A$ and not the complex conjugate transpose. Since we are using real valued basis functions in the interior domain
and therefore also for the boundary discretization, there is no difference for the surface matrices. The radiation
matrices presented later are complex valued.}

The coupling of the matrix $\Mae -\omega^2 \Mbe$ with a system matrix $\Mai -\omega^2 \Mbi$ of the interior problem becomes natural due to $\bax_j(0)=\delta_{1j}$. 
\dfo{A typical surface finite element of a $H^1$-conforming triangular finite element can be constructed such that (recall $\tilde \surf=(-\WGTH,\WGTH)$)
\begin{subequations}
\begin{align}
 \bay_1(-\WGTH)&=1,\qquad \bay_l(-\WGTH)=0,\quad l=2,\dots,N^{\rm trans},\\
 \bay_2(\WGTH)&=1,\qquad \bay_l(\WGTH)=0,\quad l=1,3,\dots,N^{\rm trans}.
\end{align}
\end{subequations}
Hence, we can associate the basis function $\bax_1 \otimes \bay_1$ to the vertex $\bV_1:=(0,-\WGTH)$ and $\bax_1 \otimes \bay_2$ to the vertex $\bV_2:=(0,\WGTH)$, which are the corners of the waveguide. All other basis functions vanish at 
these vertices. Similarly, the basis functions $\bax_1 \otimes \bay_l$ for $l=3,\dots, N^{\rm trans}$ can be associated
to the edge $\overline{\bV_1\bV_2}$, the basis functions  $\bax_j \otimes \bay_1$ for $j=2,\dots, N^{\rm long}$ to the
lower boundary of the waveguide $\WG$ and the basis functions  $\bax_j \otimes \bay_2$ for $j=2,\dots, N^{\rm long}$ to
the upper boundary of the waveguide $\WG$. The basis functions  $\bax_j \otimes \bay_l$ for $j=2,\dots, N^{\rm long}$ and $l=3,\dots,N^{\rm trans}$ are interior basis functions and vanish at $\partial \WG$.}
Therefore, we have constructed an (in)finite element for each waveguide which can be easily combined with the finite element discretization of the interior problem.
In more details and in 3d this construction is shown in \cite[Sec.~4.1]{Nannenetal:13}.

\dfo{
Recall from Sec.~\ref{sec:setting} that the waveguides are given by $\WG_\ell=\emotion_\ell\left((0,\infty)\times\tilde\surf_\ell\right)$ for $\ell=1,\dots,\nrWG$ with an Euclidean motion $\emotion_\ell$. For a solution $\bu$ to \eqref{eq:lame} we define
\begin{equation}
 \bue:=\left(\bue_1,\dots,\bue_\nrWG \right)^\top\qquad \text{with }\quad \bue_\ell:=\bu|_{\WG_\ell}\circ \emotion_\ell,\quad \ell=1,\dots,\nrWG,
 \end{equation}
and the bilinear forms  
\begin{equation}
 \Bae{\bue}{\bve}:=\sum_{\ell=1}^\nrWG \Baell{\bue_\ell}{\bve_\ell},\quad  
 \Bbe{\bue}{\bve}:=\sum_{\ell=1}^\nrWG \Bbell{\bue_\ell}{\bve_\ell}
\end{equation}
for $ \bve:=\left(\bve_1,\dots,\bve_\nrWG \right)^\top$. The bilinear forms for each waveguide are given by
\begin{subequations}
\begin{align}
  \Baell{\bue_\ell}{\bve_\ell}&:=\int_{\setR_+\times \tilde\surf_\ell} \left(2 \mu \bepsilon_\ell (\bue_\ell):\bepsilon_\ell (\bve_\ell) 
  + \lambda \dv_\ell \bue_\ell \, \dv_\ell \bve_\ell\right) \,d(\WGx,\WGy), \label{eq:BilFaWGell}\\
  \Bbell{\bue_\ell}{\bve_\ell}&:=\int_{\setR_+\times \tilde\surf_\ell} \rho \, \bue_\ell \cdot \bve_\ell \, d(\WGx,\WGy)
\end{align}
with the differential operator
\begin{equation}
 \nabla_\ell \bv_\ell := \nabla \emotion_\ell \left(\nabla \bv\right) \circ \emotion_\ell\qquad\text{for }\quad  \bv_\ell:=\bv \circ \emotion_\ell,  
\end{equation}
and $\bepsilon_\ell(\bv_\ell):=\frac{1}{2}(\nabla_\ell \bv_\ell + (\nabla_\ell \bv_\ell)^\top)$, $\dv_\ell \bv_\ell :=\nabla_\ell \cdot \bv_\ell$.
\end{subequations}
Of course, the bilinear forms in \eqref{eq:BilFaWGell} are more complicated as the one in \eqref{eq:BilFaWG} without an Euclidean motion. Nevertheless, using the same radiation matrices as in the last subsection and surface matrices with $\tilde \surf_\ell$ instead of $\tilde \surf$ and surface basis functions
\begin{equation*}
 \bayell_l(\WGy):=\left(\bai_{j_\ell(l)} \circ \emotion_\ell\right) (0,\WGy),\qquad \WGy \in \tilde \surf_\ell,\quad l=1,\dots, N^{\rm trans}_\ell,
\end{equation*}
we arrive at system matrices $\Mae_\ell -\omega^2 \Mbe_\ell \in \setC^{2 N^{\rm long}_\ell N^{\rm trans}_\ell \times 2 N^{\rm long}_\ell N^{\rm trans}_\ell}$ with $\Mae_\ell$ and $\Mbe_\ell$
similar to those in \eqref{eq:IEMatrices}. The coupling of $\Mae_\ell -\omega^2 \Mbe_\ell$ to $\Mai -\omega^2 \Mbi$ is the same as in the last section. In particular, all interface terms 
in $\Baux{\bui}{\bvi}$ cancel out with interface terms from the waveguides.}

\dfo{Formally, we can write down the variational formulation of the scattering problem \eqref{eq:lame} as follows: Find $(\bui, \bue)^\top$ 
such that
\begin{equation}\label{eq:VarForm}
\Ba{\smat{\bui\\ \bue}}{\smat{\bvi\\ \bve}} -\omega^2 \Bb{\smat{\bui\\ \bue}}{\smat{\bvi\\ \bve}}=\Lf{\bvi}
\end{equation}
for all $(\bvi, \bve)^\top$ with 
\begin{subequations}
 \begin{align}
  \label{eq:defBa} \Ba{(\bui, \bue)^\top}{(\bvi, \bve)^\top}&:=\Bai{\bui}{\bvi}+\sum_{\ell=1}^\nrWG \Baell{\bue_\ell}{\bve_\ell},\\
  \label{eq:defBb} \Bb{(\bui, \bue)^\top}{(\bvi, \bve)^\top}&:=\Bbi{\bui}{\bvi}+\sum_{\ell=1}^\nrWG \Bbell{\bue_\ell}{\bve_\ell},\\
  \Lf{\bvi} &:= \int_{\Oi} \bfc \cdot \bvi \, d\bx + \int_{\partial \Oi \setminus \left(\bigcup_{\ell=1}^L \surf_\ell \right)} \bg \cdot \bvi \, ds.   
 \end{align}
\end{subequations}
In the presence of an incoming wave $\bu^{\rm inc}$, we have to add the jumps in the Neumann data
$\sum_{\ell=1}^\nrWG \int_{\surf_\ell} \left(\sigma(\bu^{\rm inc})\cdot\bn\right)\cdot \bvi \, ds$
to the linear form and to incorporate the jumps in the Dirichlet data $\bu^{\rm inc}|_{\surf_\ell}$ to the coupling of $\bui$ and $\bue_\ell$. 
The discrete matrix form of \eqref{eq:VarForm} is given by
\begin{equation}
 \left(\Ma - \omega^2 \Mb\right) {\bf s}= {\bf F}
\end{equation}
with the solution vector ${\bf s}\in \setC^{2N}$ with $N=N^{\rm int}+\sum_{\ell=1}^\nrWG (N^{\rm long}_\ell-1) N^{\rm trans}_\ell$, the matrices 
$\Ma,\Mb \in \setC^{2N \times 2N}$ constructed out of $\Mai$, $\Mbi$, $\Mae_\ell$, and  $\Mbe_\ell$ and a right hand side ${\bf F} \in \setC^{2 N}$.}

\section{Resonance problems}\label{sec:resprob}
Up to now, we have only considered diffraction problems.
We now focus on the resonance problem \dfn{\eqref{eq:res_classic}}.
\dfn{
Since the modal radiation condition from Sec.~\ref{Sec:ModalAnalysis}
is (more ore less well) defined only for positive frequencies $\omega$,
so is the resonance problem~\eqref{eq:res_classic}.
A meaningful extension to complex valued frequencies consists
in a holomorphic extension of~\eqref{eq:res_classic}.
If $\setWN(\omega)\subset\curve_{\pol,\poli}^+$ holds true for all
$\omega$ in a real interval $(a,b)$, then the pole condition \eqref{eq:pc}
with respect to the parameters $\pol,\poli$ is equivalent to the modal
radiation condition on this interval. Since the pole condition is independent
of the frequency, it is straightforward to holomorphically extend the
resonance problem from the interval $(a,b)$ into the complex plane through
the pole condition.}

\dfn{However, the dependence of this extension on $\pol,\poli$ is
not obvious. The crucial property to investigate is the Fredholmness of the
holomorphic extension. Due to similar results for acoustic waveguides in \cite{HohageNannen:15}, we conjecture that the essential spectrum,
i.e.~the set of frequencies $\omega$ for which~\eqref{eq:res_classic}
equipped with the pole condition is not Fredholm, takes the form
\begin{equation}
\sigma_{\rm ess}(\pol,\poli):=\left\{\omega\in\setC\colon \iu\kappa_n(\omega)\in\curve_{\pol,\poli}\text{ for all  }n\in\setN\right\}.
\end{equation}
For complex frequencies $\omega$ the wavenumbers $\kappa_n(\omega)$ are thereby
defined in the same way as for real valued frequencies, i.e.
as eigenvalues of a quadratic eigenvalue problem (see Rem.~\ref{Rem:TheoModes}).
The set $\sigma_{\rm ess}(\pol,\poli)$ can be interpreted as the union of $\setN$ curves,
i.e.~$\sigma_{\rm ess}(\pol,\poli)=\bigcup_{n\in\setN}\{\omega\in\setC\colon\iu\kappa_n(\omega)\in\Gamma_{\pol,\poli}\}$.
Due to the point symmetry of $\Gamma_{\pol,\poli}$ and of the set of wavenumbers with respect to $0$,
each curve starts at a frequency $\omega$, such that there exists a vanishing wavenumber $\kappa_n(\omega)=0$.}


\begin{figure}
\centering
\subfigure[\label{fig:ess_spec1}$\pol=\poli$ as in Fig.~\ref{fig:ComplexWNrsHP}]{\resizebox{0.49\textwidth}{!}{\includegraphics{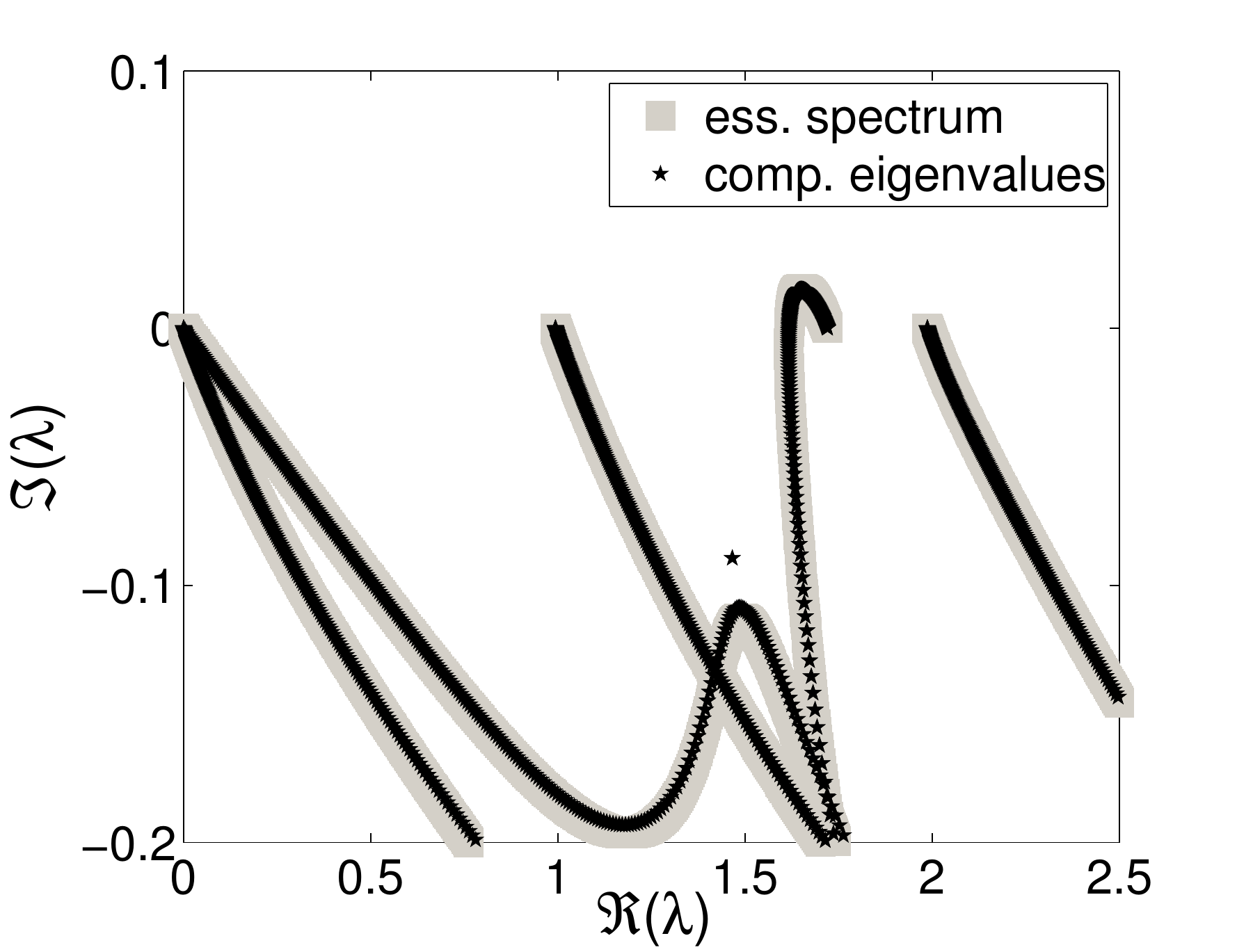}}} \hfill
\subfigure[\label{fig:ess_spec2}$\pol, \poli$ as in Fig.~\ref{fig:ComplexWNrsCurved}]{\resizebox{0.49\textwidth}{!}{\includegraphics{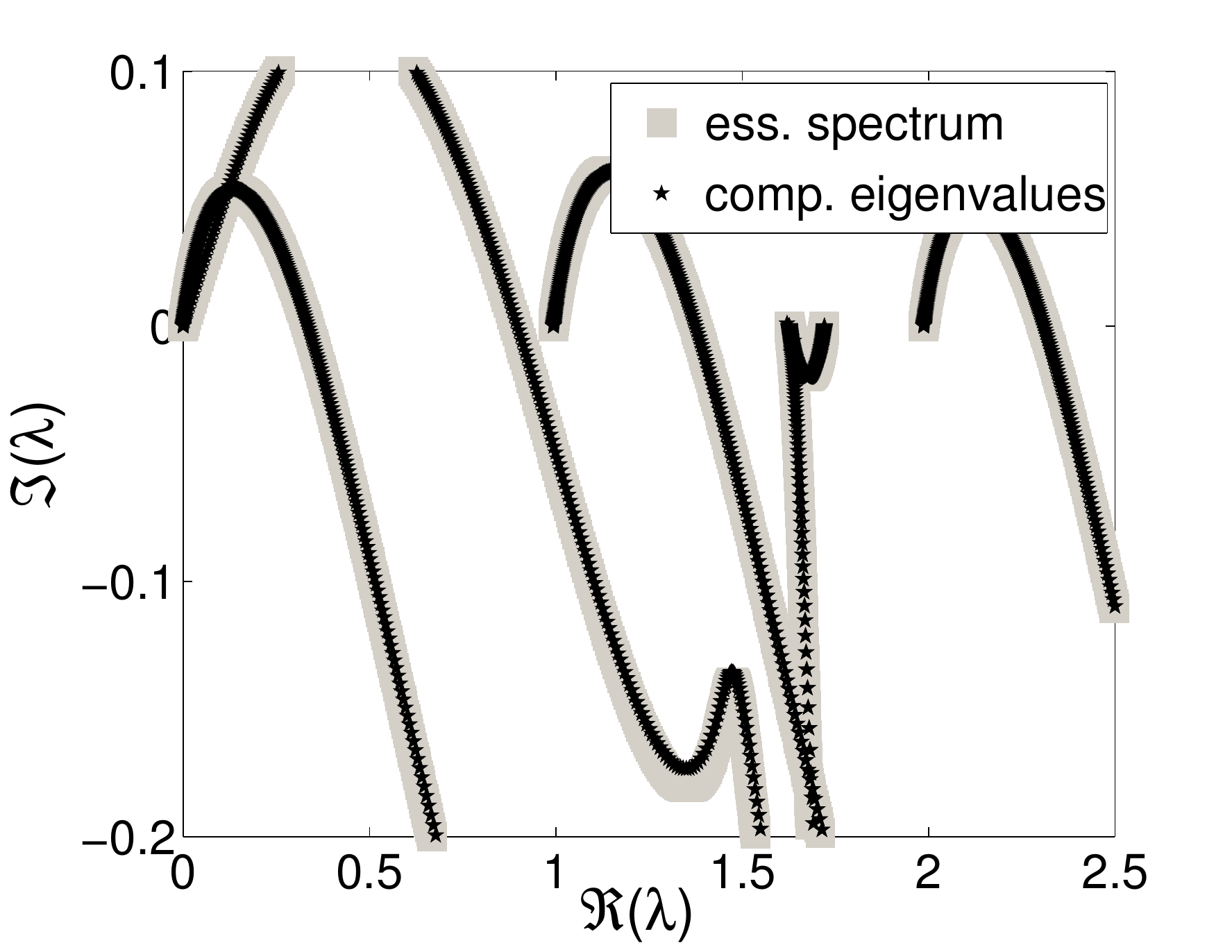}}}
\caption{Eigenvalues of \dfn{\eqref{eq:res_classic} for $\Omega=\setR \times (-1,1)$, $\rho=E=1$ and $\nu=0.25$ compared with the essential
spectrum $\sigma_{\rm ess}(\pol,\poli)$.} \dfo{A discretization of an essential spectrum, depending on $\pol,\poli$, can be observed.}}
\label{fig:ess_spec}
\end{figure}

\dfn{
To support these conjectures, we give a numerical example in Fig.~\ref{fig:ess_spec}
for a waveguide $\setR \times (-1,1)$ with $\rho=E=1$ and $\nu=0.25$.
For a given $\kappa \in -\iu \curve_{\pol,\poli}$ the corresponding $\omega$
for which wavenumbers $\kappa_n(\omega)\in-\iu\curve_{\pol,\poli}$ exist,
are the square roots of a linear eigenvalue problem given in \cite{Pagneux}.
Moreover, $-\iu \curve_{\pol,\poli}$ can be explicitly parameterized with
$-\iu \gamma_{\pol,\poli}$ due to \eqref{eq:ParamCurve}. The set $\sigma_{\rm ess}(\pol,\poli)$
can therefore be numerically computed with a standard eigenvalue solver.
For our example we chose 1000 sample points $\kappa=-\iu\gamma_{\pol,\poli}(r)\in\iu\curve_{\pol,\poli}$
with equidistant $r\in(-5,0)$ to compute $\sigma_{\rm ess}(\pol,\poli)$.
In Fig.~\ref{fig:ess_spec} we compare it with the computed eigenvalues of \eqref{eq:res_classic},
where we used a high order finite element method for $\Oi=(0,5)\times (-1,1)$
with a triangular mesh with maximal meshwidth $h=0.1$ and polynomial order $p=12$.
The two waveguides $(5,\infty)\times (-1,1)$ and $(-\infty,0)\times (-1,1)$ are
discretized using Hardy space infinite elements with $N^{\rm long}=200$.
For the generalized linear eigenvalue problem we have used a shift-and-invert
Arnoldi algorithm with four different shifts and a Krylov space with dimension $2000$.
}

\dfn{
It can be seen in Fig.~\ref{fig:ess_spec} that almost all of the computed
eigenvalues $\lambda \in \setC$ fit very well to $\sigma_{\rm ess}(\pol,\poli)$.
Thus their interpretation as the discretization of an essential spectrum is
reasonable. The isolated point in Fig.~\ref{fig:ess_spec1} has to be interpreted
as a resonance. We clearly see in Fig.~\ref{fig:ess_spec} that different parameter
choices $\pol,\poli$ lead to different essential spectra, as predicted through
$\sigma_{\rm ess}(\pol,\poli)$. We also observe a dependence of the discrete spectrum
on the parameters. This is because the different parameters correspond to different
choices of branches of the solution operator. This behavior will be discussed
in more detail in Sec.~\ref{sec:MP3}.}

\dfn{
Moreover, we observe in Fig.~\ref{fig:ess_spec} that some curves of the essential spectrum
are partially located in the quadrant $\{z\in\setC\colon\Re(z),\Im(z)>0\}$.
Each such curve admits in addition to its starting point a supplementary intersection
$\omega^*$ with the real axis. This means that there exists a wavenumber $\kappa_n(\omega)$
crossing over $-\iu \curve_{\pol,\poli}$ as $\omega$ moves along the real interval
$(\omega^*-\epsilon,\omega^*+\epsilon)$ (for some $\epsilon>0$).
Since the pole condition selects wavenumbers in $-\iu \curve_{\pol,\poli}^+$, it
can't be equivalent to the modal radiation condition for
both intervals $(\omega^*-\epsilon,\omega^*)$, $(\omega^*,\omega^*+\epsilon)$
simultaneously. A inspection of such cases shows, that for intervals $I$ with
$\sigma_{\rm ess}(\pol,\poli)\cap\{x+\iu y\colon x\in I,y>0\}\neq\emptyset$ the pole
condition is never equivalent to the modal radiation condition. 
Note, that such an interval exists in Fig.~\ref{fig:ess_spec1} indicating the presence of a backward propagating mode for $\omega \in (1.58,1.65)$ (compare with Fig.~\ref{fig:disprelfrequencyrange}). 
Since the Hardy space infinite element method only depends on $\pol$ and $\poli$, for this finding no calculation of wavenumbers or waveguide modes was needed.
}

\dfn{
It is worth mentioning that the discrete and essential spectra can also be
interpreted in terms of the solution operator to \eqref{eq:res_classic} with
respect to the pole condition. Resonances correspond to poles and the essential
spectrum corresponds to branch cuts of the solution operator.}

\dfn{We conclude, that the computed eigenvalues of \eqref{eq:res_classic}
have to be interpreted carefully. They might be part of the discretization
of an essential spectrum. If this is not the case, they are approximations
to the resonances of \eqref{eq:res_classic}. These resonances depend on the
parameter choice and belong to different Riemann sheets.
The relevance of resonances for scattering problems will be
investigated further in Sec.~\ref{sec:MP2}}
\dfn{Note, that the discretization of \eqref{eq:res_classic} with Hardy
space infinite elements leads to a generalized linear eigenvalue problem,
which can be solved using a standard linear eigenvalue solver. This is
one of the advantages of the presented method over classical modal methods.}

\section{Numerical results}\label{sec:numres}
Since rigorous convergence results are not available for Hardy space elements in the context of elastic waveguide problems, we report here on numerical experiments. The convergence results 
for acoustic waveguide problems in \cite{HohageNannen:15} and of the model problem in \cite{HallaHohageNannenSchoeberl:14} indicate a super-algebraic convergence 
rate with respect to the number of unknowns in the infinite direction. The first numerical experiment for a diffraction problem with a known solution supports this conjecture. 
The second numerical experiment underlines the relation between a scattering problem and the corresponding resonance problem. 
The last experiment illustrates the dependency of the computed resonances on the method parameters.

The computations were all made with the finite element package netgen/NGSolve \cite{netgen,ngsolve:14} and the module ngs-waves \cite{ngs-waves} containing the source code for Hardy space methods. \dfn{For the resonance problems we have used a standard shift-and-invert Arnoldi method using MUMPS or PARDISO as direct solver.}

\subsection{Convergence test}
\begin{figure}[t]
\centering
\subfigure[\label{fig:error1w}]{\resizebox{0.45\textwidth}{!}{\includegraphics{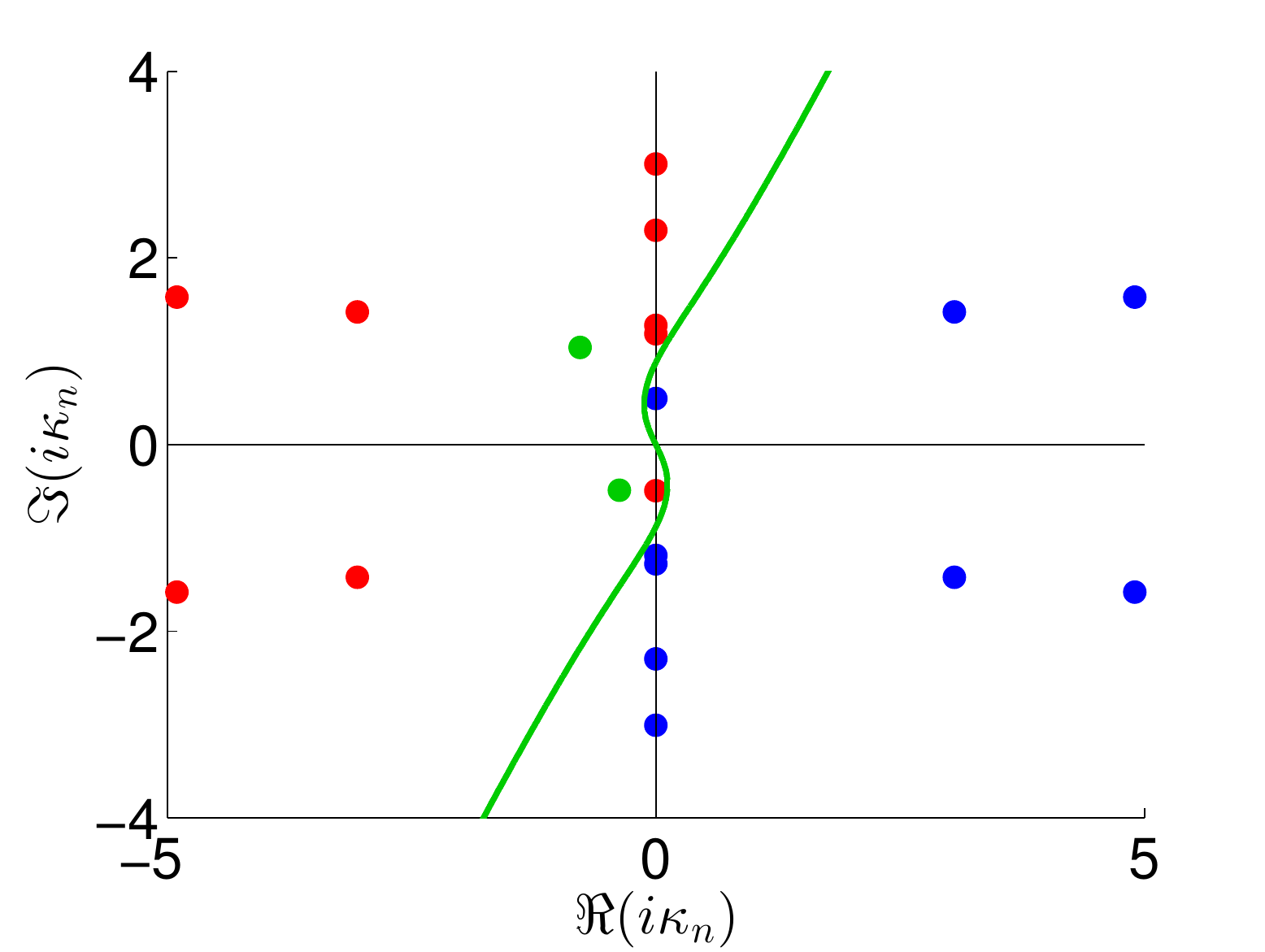}}} 
\subfigure[\label{fig:error1a}]{\resizebox{0.45\textwidth}{!}{\includegraphics{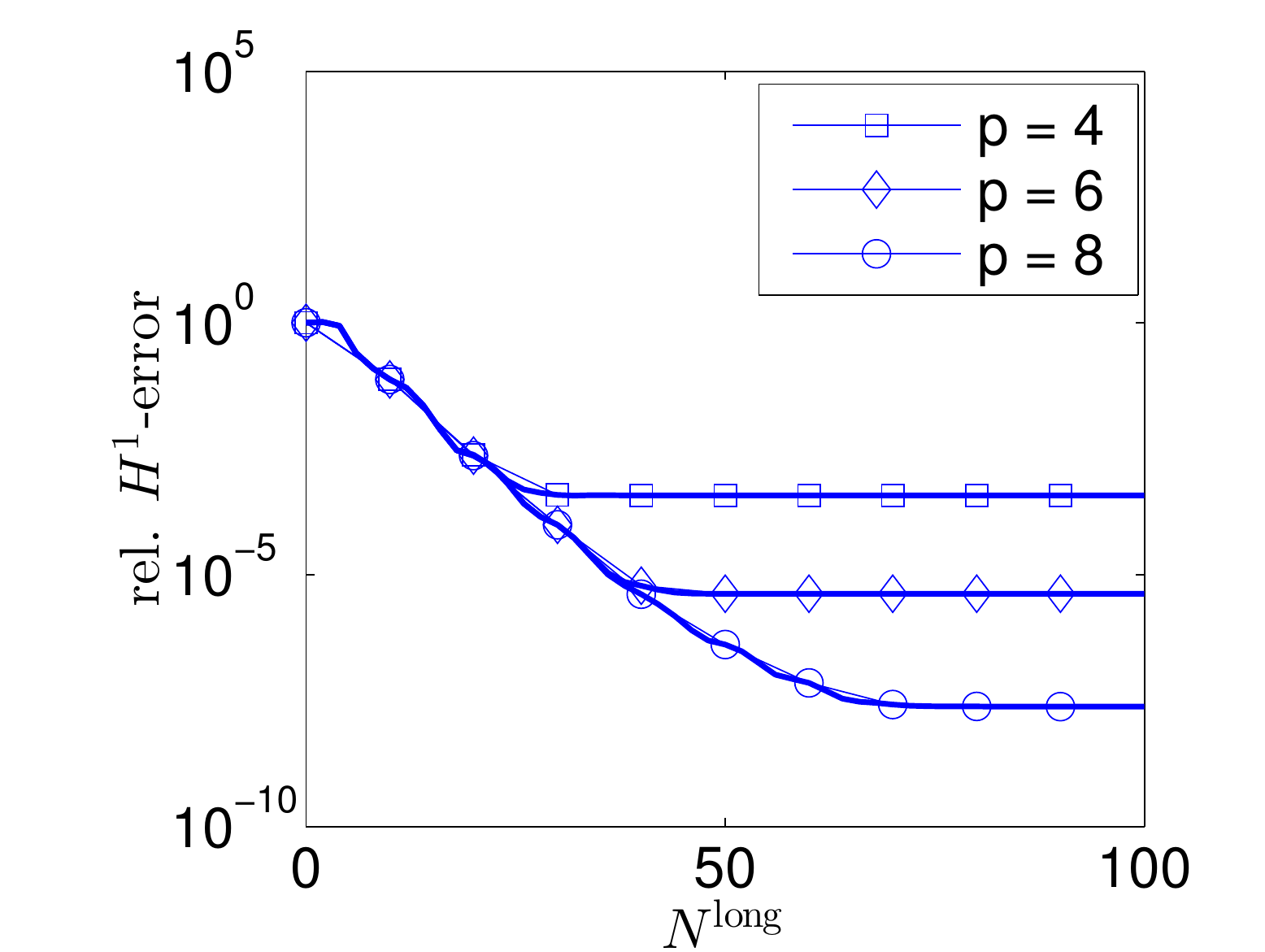}}} 
\caption{Left panel: outgoing (red) and incoming (blue) $\iu\kappa_n$ for $\omega=1.66$ \dfo{(squares) and $\omega=1.7205$ (diamonds)}. Right panel: relative $\left(H^1(\Oi)\right)^2$-error of $\bu-\bu_{\rm ref}$ w.r.t. the number of unknowns $N^{\rm long}$ in the infinite direction for different polynomial orders.}
\label{fig:error1}
\end{figure}

\begin{figure}[!b]
\centering
\subfigure[\label{fig:error2a} $\omega \in (1.57, 1.78)$, $p=5$]{\resizebox{0.45\textwidth}{!}{\includegraphics{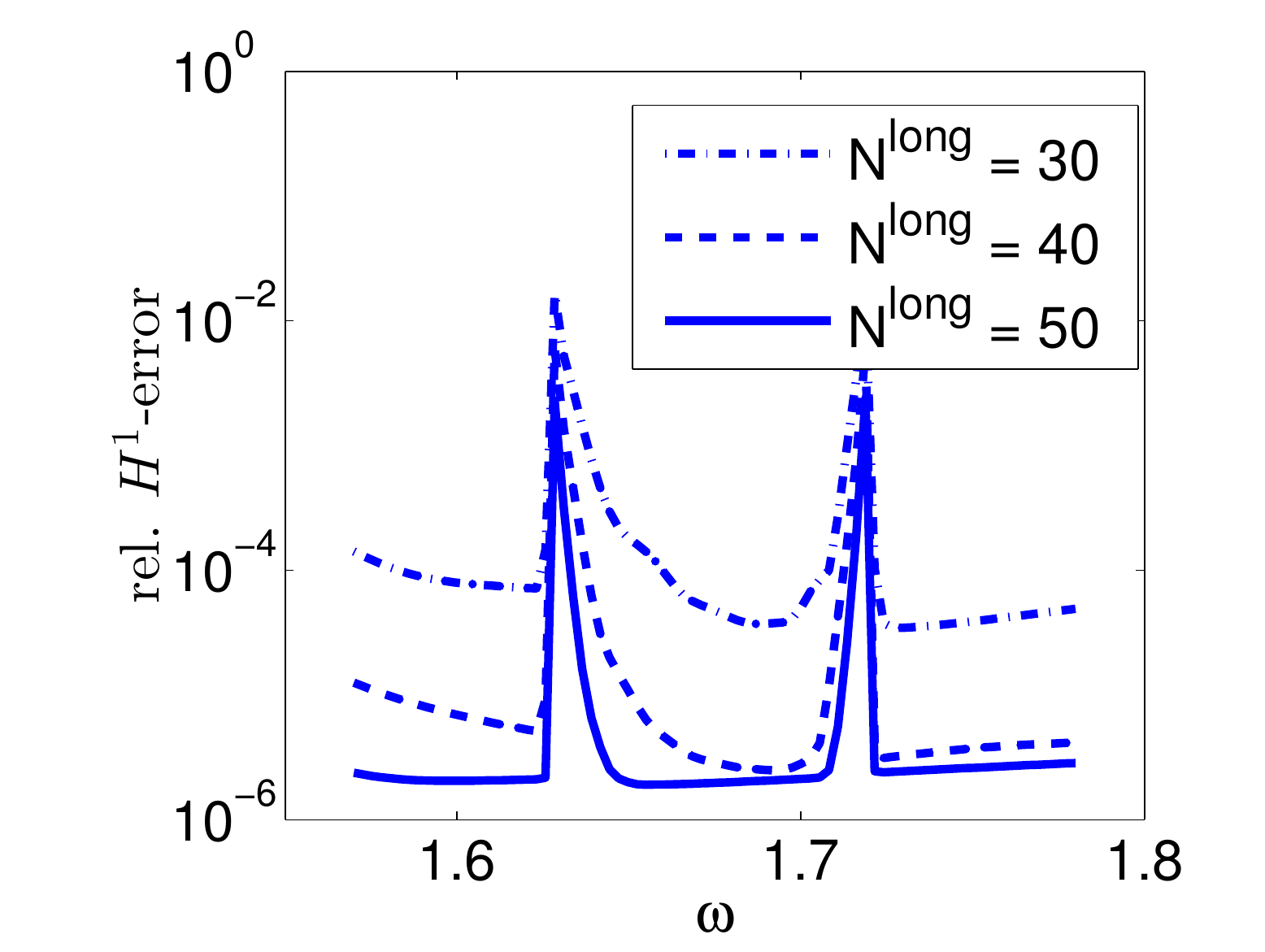}}}
\subfigure[\label{fig:error2b} error w.r.t. to $N^{\rm long}$ for $p=10$]{\resizebox{0.45\textwidth}{!}{\includegraphics{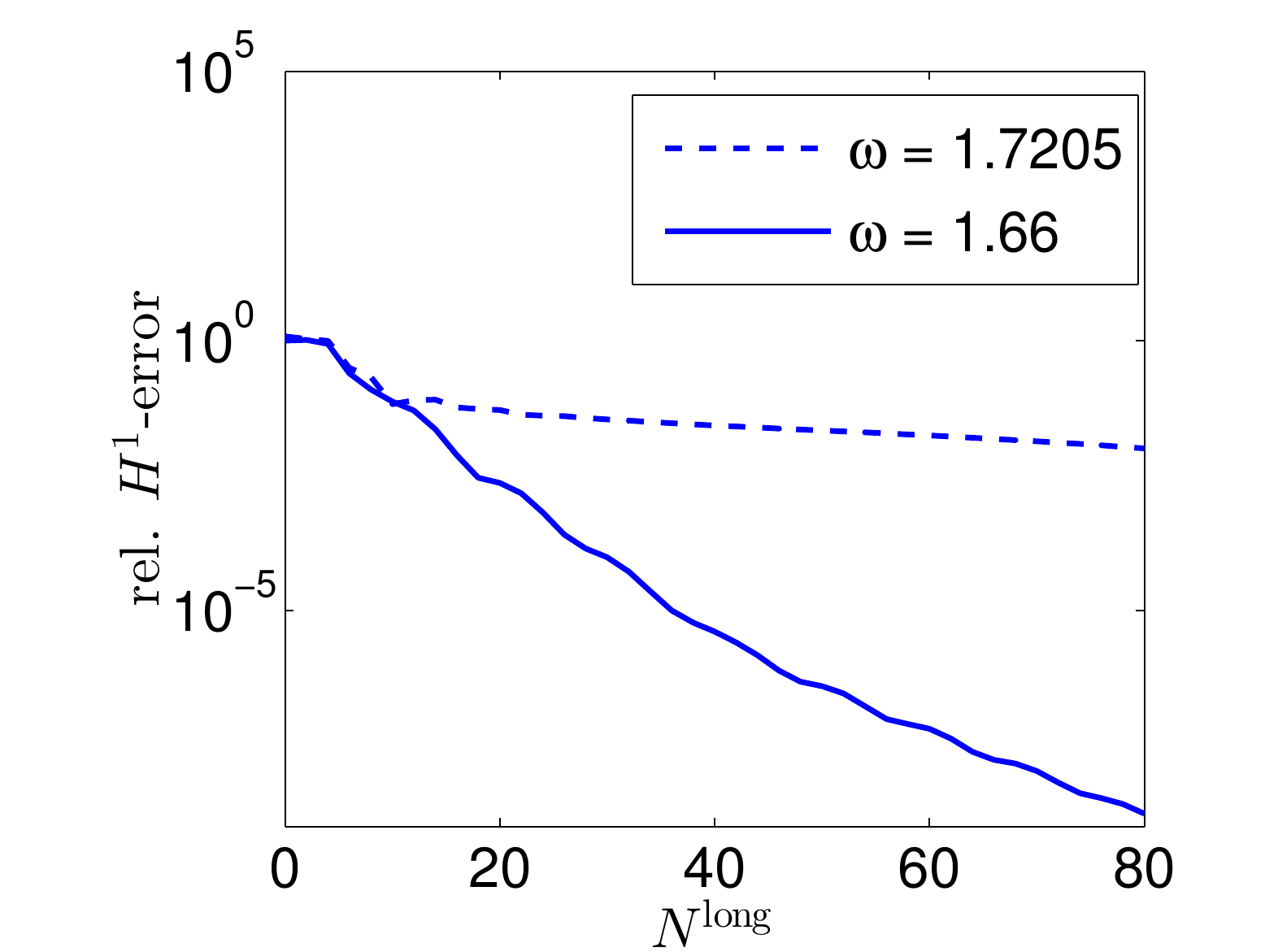}}}
\caption{relative $\left(H^1(\Oi)\right)^2$-error of $\bu-\bu_{\rm ref}$ for different frequencies $\omega$}
\label{fig:error2}
\end{figure}
For unperturbed waveguides $W=\setR^+\times(-\WGTH,\WGTH)$ the Rayleigh-Lamb modes are known to be (semi-)explicit solutions to \eqref{eq:modal} (see e.g. \cite{Achenbach73}): 
We define the longitudinal wave speed by $c_L:=\sqrt{(\lambda+2\mu)/\rho}$, the transversal wave speed by $c_T:=\sqrt{\mu/\rho}$), and for $\kappa\in \setC$  the dispersion relations 
 \begin{subequations}\label{eq:disprel}
\begin{align}
\label{eq:disprel0} F_{S}(\kappa) &:= 4\kappa^2\alpha \beta \sin (\alpha \WGTH) \cos (\beta \WGTH) + (\kappa^2-\beta ^2)^2 \cos (\alpha \WGTH) \sin (\beta \WGTH),  \\
\label{eq:disprel1} F_{A}(\kappa) &:= 4\kappa^2\alpha \beta \cos (\alpha \WGTH) \sin (\beta \WGTH) + (\kappa^2-\beta ^2)^2 \sin (\alpha \WGTH) \cos (\beta \WGTH),
\end{align}
\end{subequations}
with $\alpha:=\sqrt{\omega^2/c_L^2-\kappa^2}$ and $\beta:=\sqrt{\omega^2/c_T^2-\kappa^2}$. If $\kappa$ is a root of $F_{S}(\kappa)=0$ or $F_{A}(\kappa)=0$, then 
\begin{equation*}
 \bu^S(\WGx,\WGy):=e^{i \kappa \WGx} \bw^S(\WGy), \quad \bu^A(\WGx,\WGy):=e^{i \kappa \WGx} \bw^A(\WGy),\qquad (\WGx,\WGy) \in W,
\end{equation*}
with
\begin{subequations}\label{eq:lambmodes}
\begin{align}
\bw^{S}(\WGy)&:=\smat{i \kappa(\kappa^2-\beta^2)^2\sin(\beta \WGTH) \cos(\alpha \WGy) + \beta 2i\kappa\alpha\sin (\alpha \WGTH)\cos(\beta \WGy) \\ 
-\alpha (\kappa^2-\beta^2)^2 \sin(\beta \WGTH)\sin(\alpha \WGy) + 2\kappa^2\alpha\sin(\alpha \WGTH)\sin(\beta \WGy)},\\
\bw^{A}(\WGy)&:= \smat{i\kappa(\kappa^2-\beta^2)^2\cos(\beta \WGTH) \sin(\alpha \WGy) - \beta 2i\kappa\alpha\cos(\alpha \WGTH)\sin(\beta \WGy) \\ 
\alpha (\kappa^2-\beta^2)^2\cos(\beta \WGTH)\cos(\alpha \WGy) + 2\kappa^2\alpha\cos(\alpha \WGTH)\cos(\beta \WGy)},
\end{align}
\end{subequations}
are the symmetric and anti-symmetric Rayleigh-Lamb modes.
 
We chose for $\WGTH=1$, Young's modulus $E=1$, density $\rho=1$ and Poisson's ratio $\nu=0.25$ the reference function
\begin{equation}
 \bu_{\rm ref}:=\sum_{j=1}^5 \bu^S_j/\|\bw^S_j\|_{L^2(\surf)} + \sum_{j=1}^4 \bu^A_j/\|\bw^A_j\|_{L^2(\surf)}
\end{equation}
with the first five symmetric and first four anti-symmetric Rayleigh-Lamb modes. The domain is given by 
$\Omega=\Oi\cup\surf\cup W$ with $\Oi=(0,15)\times(-1,1)$ triangulated with maximal mesh size $h=0.25$, $\surf=\{15\}\times(-1,1)$, and $W=(15,\infty)\times(-1,1)$.
\eqref{eq:modal} was completed with the Dirichlet boundary condition $\bu(0,\bullet)=\bu_{\rm ref}(0,\bullet)$ and the pole condition \dfn{\eqref{eq:pc}}\dfo{Def.~\ref{def:PoleCondition}}  for $\curve_{\pol,\poli}$ defined   
in Sec.~\ref{Sec:defCurve} with $\pol=-0.374158-0.488609\iu$ and $\poli=-0.775234+1.03962\iu$.

First we pick a fixed frequency $\omega=1.66$, such that there exists an outgoing wavenumber $\kappa(\omega)<0$ (see Fig.~\ref{fig:error1w}) and vary the polynomial order 
of the finite element method for $\Oi$ and the number of unknowns $N^{\rm long}$ in the infinite direction (see Fig.~\ref{fig:error1a}). \dfn{Super-algebraic}\dfo{Exponential} convergence in $N^{\rm long}$ can be observed \dfn{until the error of the finite element discretization of the interior problem, which depends on the polynomial order, is reached}.

In Fig.~\ref{fig:error2a} we fixed the uniform polynomial degree to $p=5$ and varied $\omega \in (1.57, 1.78)$. At the 
frequencies $\omega\approx 1.6260,\omega \approx 1.7206$ the method fails, as already mentioned in Sec.~\ref{Sec:CondCurve}, since for these frequencies one outgoing wavenumber 
coincides with an incoming one. Apart from these frequencies, the relative error in $\Oi$ is small. But if we chose 
$\omega=1.7205$ such that there exists one wavenumber $\kappa_n\approx 0$ \dfo{(see the diamonds in \ref{fig:error1w} near to the origin)}, the convergence rate w.r.t. $N^{\rm long}$ 
is very poor (see \ref{fig:error2b} with $p=10$). Techniques to overcome this problem were developed for acoustic waveguide problems in \cite{HohageNannen:15}. The application to elastic 
waveguide problems is intended for future research.

\subsection{Cavity resonances}
\label{sec:MP2}
\begin{figure}[b]
\centering
\subfigure[\label{fig:geom2} triangulation of $\Oi$]{\resizebox{0.45\textwidth}{!}{\includegraphics{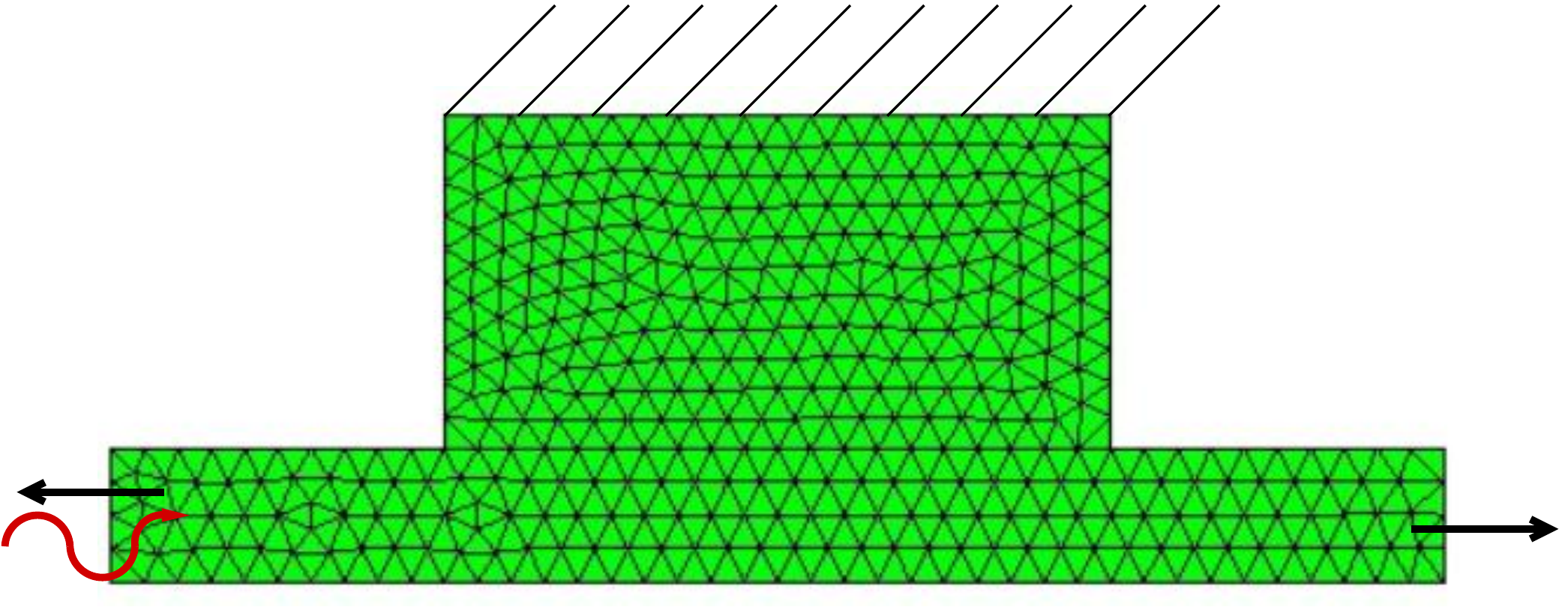}}}
\subfigure[\label{fig:stress} stress in $\Omega_1$ w.r.t. $\omega$]{\resizebox{0.3\textwidth}{!}{\includegraphics{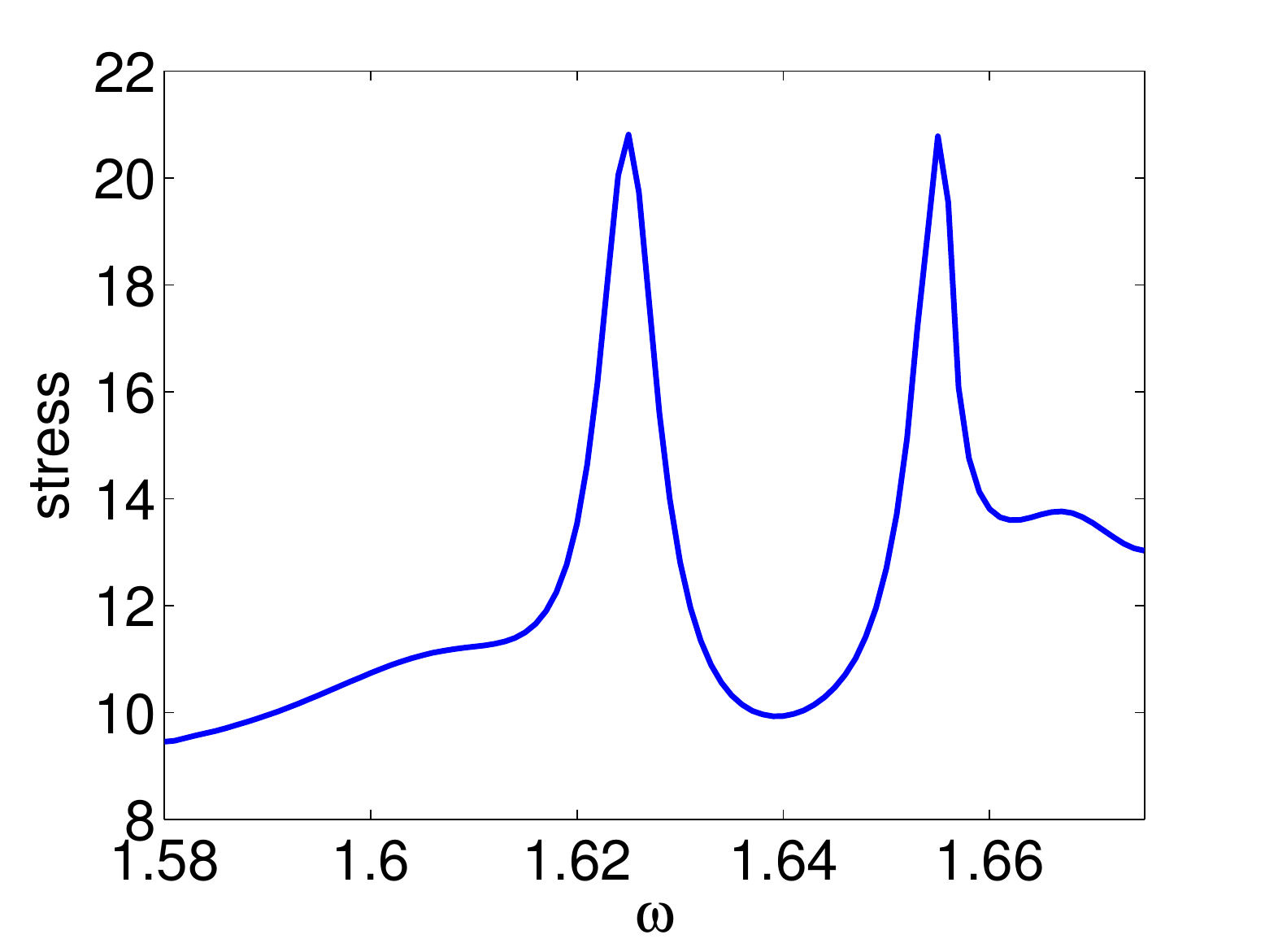}}}
\caption{scattering of a wave signal by a cavity with Dirichlet boundary conditions}
\end{figure}

We consider a waveguide $\setR\times (-1,1)$ with a cavity $\Omega_1:=(-5,5)\times [1,6)$
attached to some wall at the homogeneous Dirichlet-boundary $\partial\Omega_D=(-5,5)\times\{6\}$.
On all other parts of the boundary we assume homogeneous Neumann boundary conditions.
In Fig.~\ref{fig:geom2} the triangulation of the interior domain $\Oi=(-10,10)\times(-1,1)\cup \Omega_1$ is given. 

\begin{figure}[!b]
\centering
\subfigure{\resizebox{0.35\textwidth}{!}{\includegraphics{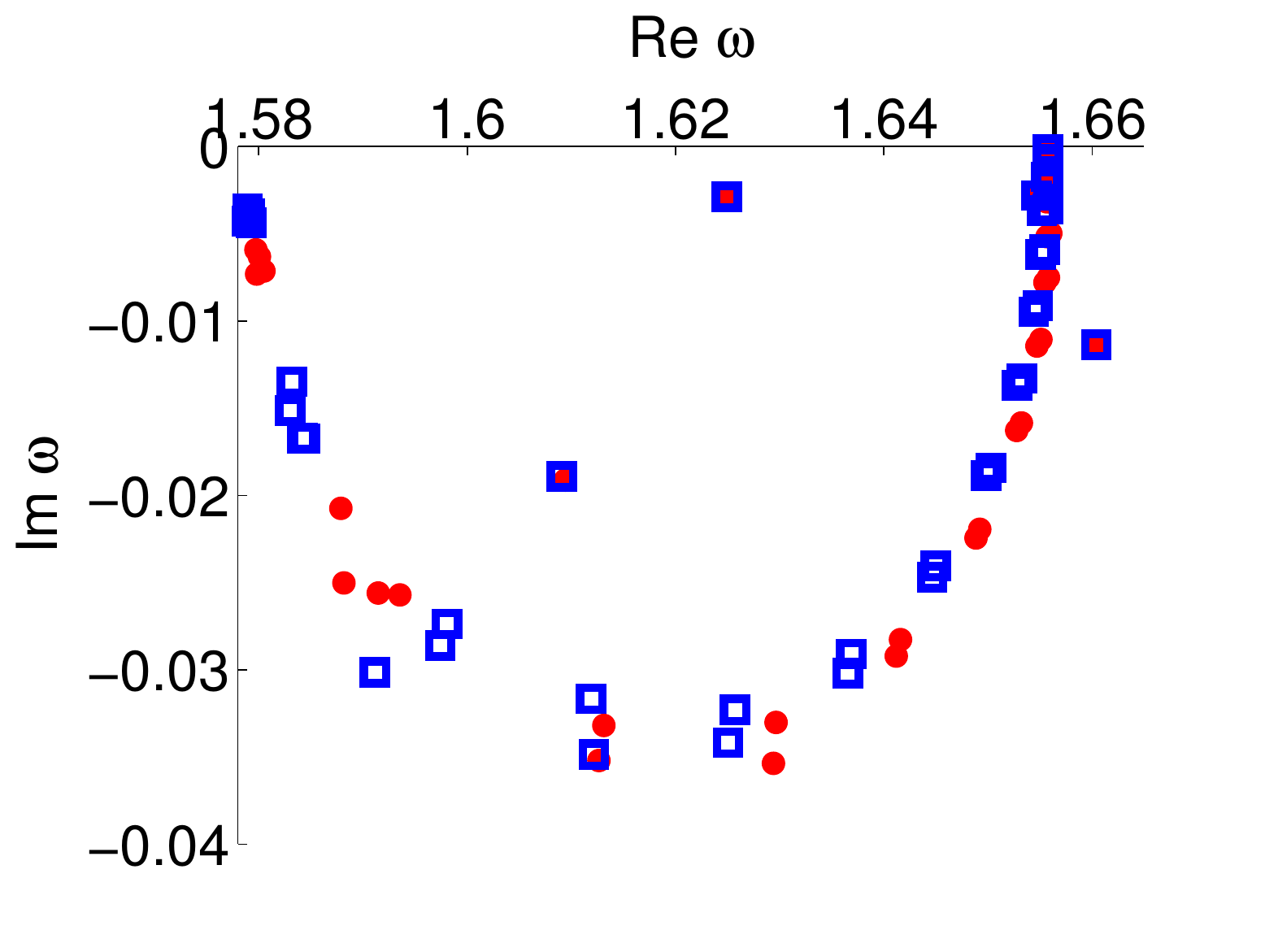}}}\qquad
\subfigure{\resizebox{0.35\textwidth}{!}{\includegraphics{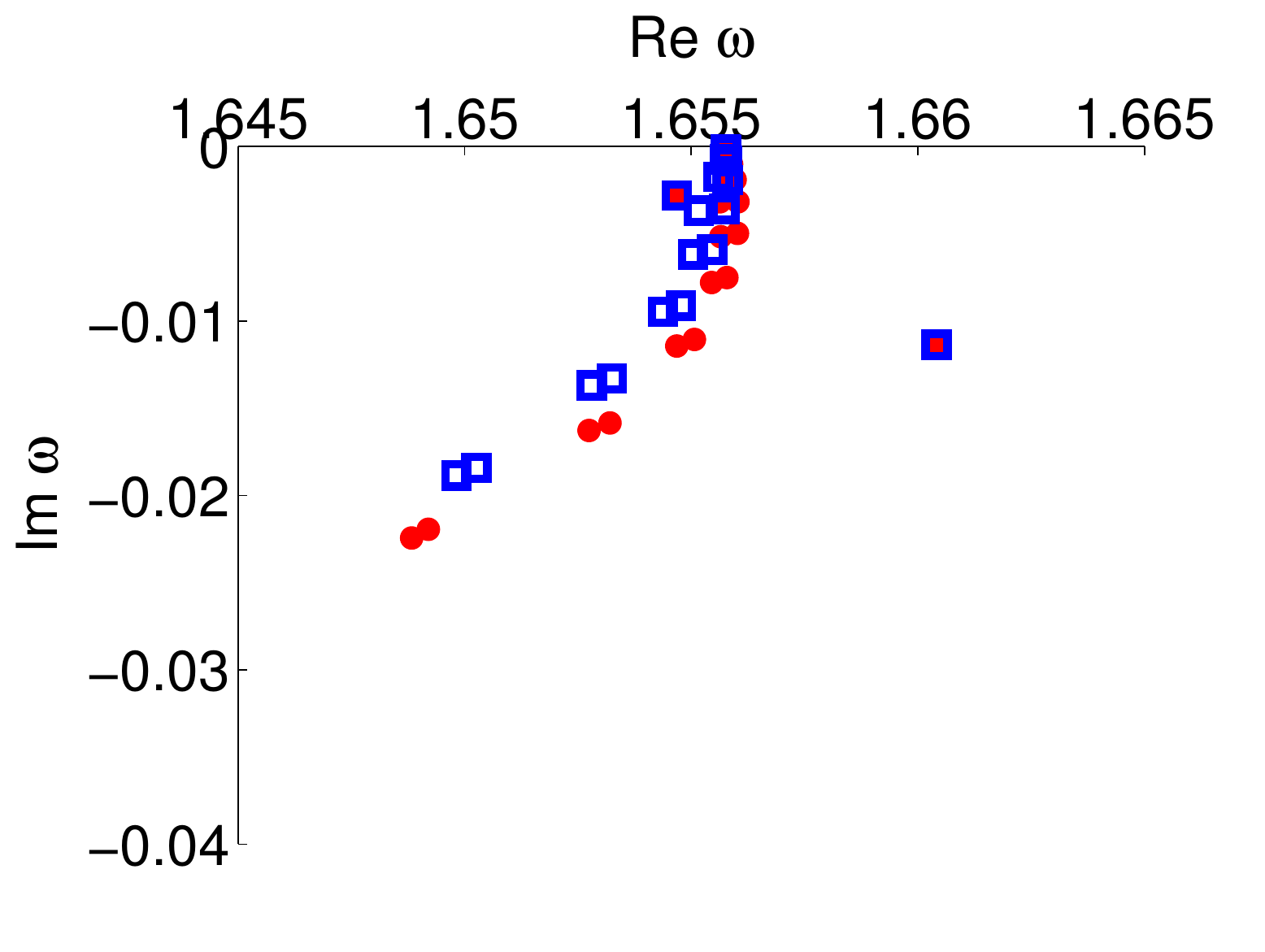}}}
\caption{The computed spectrum for \dfn{a pair or parameters $\pol,\poli$ indicated as red points,
and for a slightely perturbed pair $\tilde\pol,\tilde\poli$ indicated as blue squares.}
\dfo{two slightly different pairs of parameters $\pol,\poli$.}
The right plot is a magnification of one part of the left plot.}
\label{fig:res-prob}
\end{figure}

\dfn{First}\dfo{Second}, we have solved \dfo{the} scattering problem \eqref{eq:lame} for different frequencies
\dfo{$\omega$}\dfn{$\omega\in (1.58,1.65)$} with an incoming wave consisting of propagating Lamb modes and measured the stress
in the cavity $\Omega_1$\dfn{, which is plotted in Fig.~\ref{fig:stress}}. \dfn{We have chosen the material parameters as in the last subsection,
in order to ensure the existence of backward propagating modes in the stated frequency interval.}

\dfn{Second}\dfo{First}, we solve the resonance problem \eqref{eq:res_classic} \dfn{for the same material parameters} \dfo{for $E=\rho=1, \nu=0.2$}
\dfn{and look for resonances near the interval $(1.58,1.65)$}.
In Fig.~\ref{fig:res-prob} the computed eigenvalues for two different pairs of $\pol$ and $\poli$ are given.
\dfn{We can detect in Fig.~\ref{fig:res-prob} four eigenvalues, are separated from the remaining
spectrum. They are $\omega_1=1.609-0.019\iu$, $\omega_2=1.625-0.003\iu$, $\omega_3=1.655-0.003\iu$, and $\omega_4=1.66-0.011\iu$. Moreover these eigenvalues coincide for both calculations obtained with different method parameters. Since the essential spectrum depends on the method parameters, we can identify them as resonances.}

\dfn{If we compare Fig.~\ref{fig:stress} with Fig.~\ref{fig:res-prob}, we observe that}
the two peaks in Fig.~\ref{fig:stress} fit very well to the two computed resonances
$\omega_{\rm res}$ from Fig.~\ref{fig:res-prob} with $|\Im (\omega_{\rm res})|$ most low\dfn{, i.e.~$\omega_2$ and $\omega_3$}.
\dfn{This is not surprising, since apart from the essential spectrum we expect the solution operator to be meromorphic with respect
to the frequency with resonances being the poles.
Hence, in the neighborhood of a resonance the scattering problem is almost singular and thus
very sensitive to external forces.}

\dfo{the solution will contain backward propagating modes.}
\dfo{Since the discretization of the essential spectrum depends on $\pol$ and $\poli$,
resonances appearing in both computations are not part of an essential spectrum.
Hence, we can identity in  Fig.~\ref{fig:res-prob} four different resonances.}

\subsection{Dependency of resonances with respect to $\pol,\poli$}
\label{sec:MP3}

\begin{figure}[tb]
\centering
\includegraphics[width=0.8\textwidth]{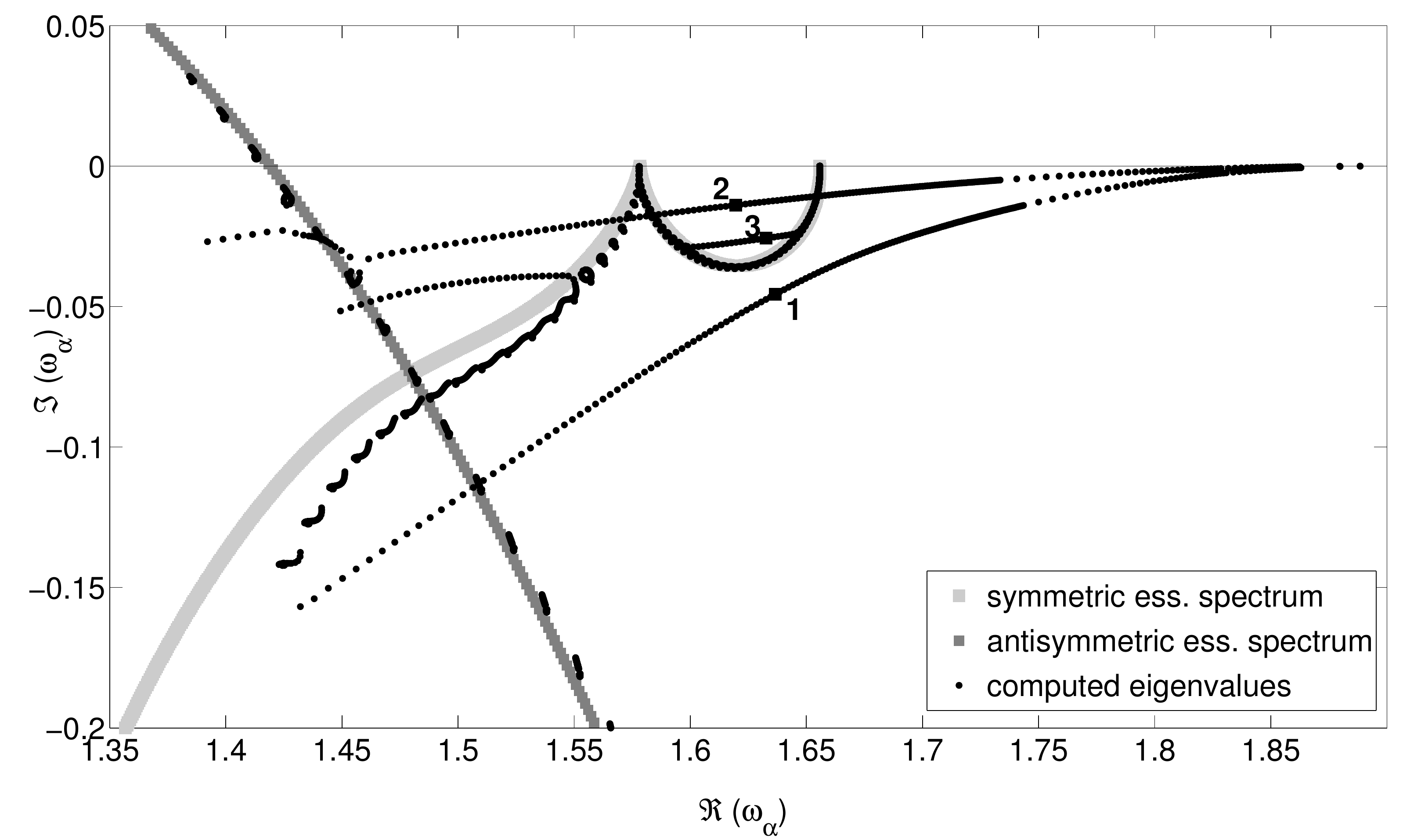}
\caption{the computed spectrum for $\pol=-0.182086-0.237784\iu$, $\poli=-1.71492+2.29978\iu$ and \dfn{varying} $\alpha\in(5,10^{10})$;
\dfn{eigenvalues near to}\dfo{ resonances in} shaded domains are discretizations of \dfn{two curves $\{\omega\in\setC\colon\iu\kappa_n(\omega)\in\Gamma_{\pol,\poli}\}$ of} the essential spectrum \dfn{(see Sec.~\ref{sec:resprob}) with $\kappa_n$ being the wavenumbers to one antisymmetric Lamb mode and to one symmetric, backward propagating Lamb mode. The squares indicate the resonances for $\alpha=11$.}}
\label{fig:res_alpha_variation1}
\end{figure}

\begin{figure}
\centering
\subfigure[\label{fig:res1comp} resonance 1 in Fig.~\ref{fig:res_alpha_variation1} ($\omega\approx 1.636 - 0.045\iu$)]{\resizebox{0.45\textwidth}{!}{\includegraphics{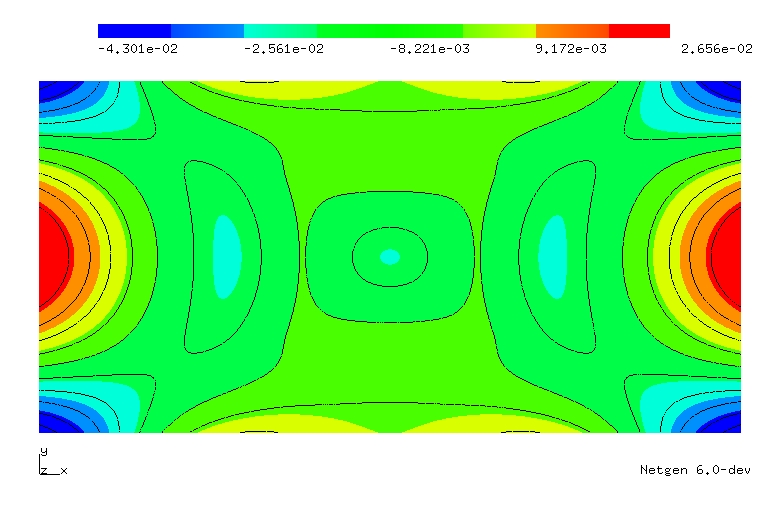}}
\resizebox{0.45\textwidth}{!}{\includegraphics{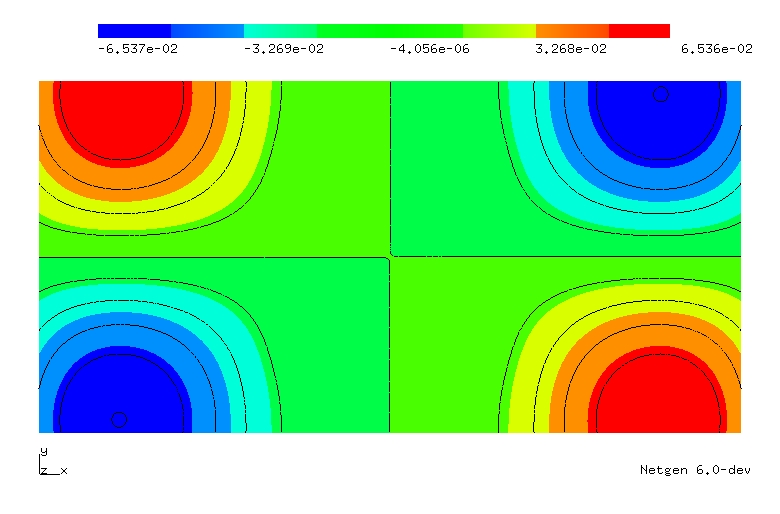}}}
\subfigure[\label{fig:res2comp} resonance 2 in Fig.~\ref{fig:res_alpha_variation1} ($\omega\approx 1.620 - 0.014\iu$)]{\resizebox{0.45\textwidth}{!}{\includegraphics{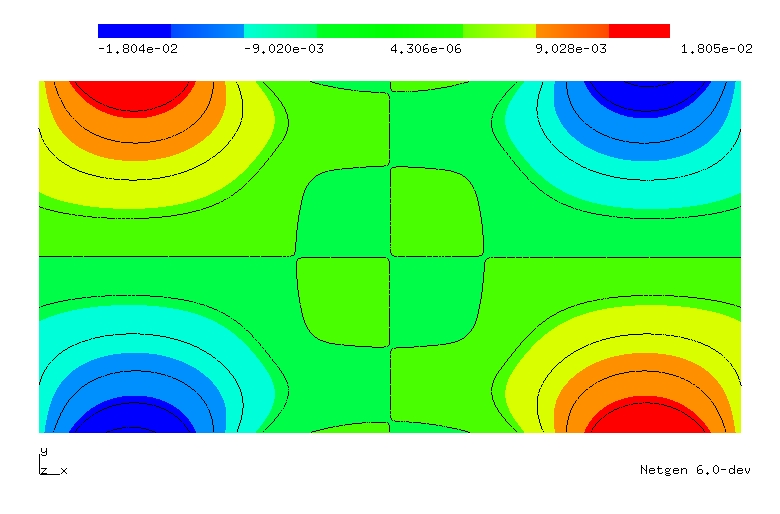}}
\resizebox{0.45\textwidth}{!}{\includegraphics{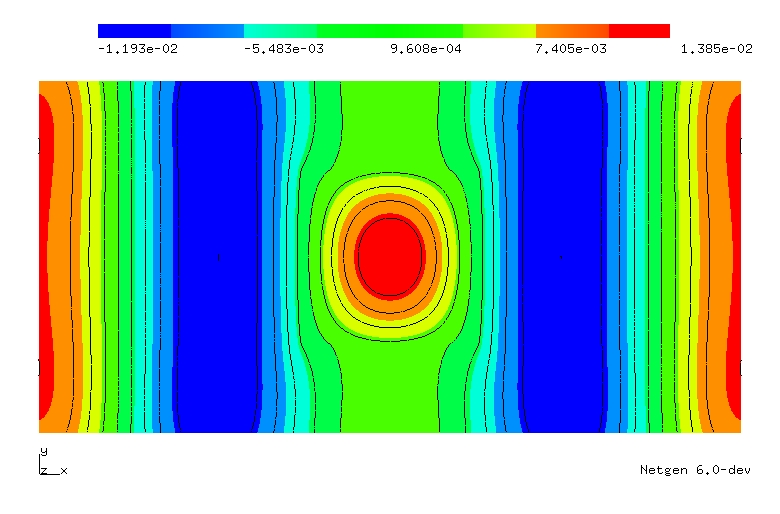}}}
\subfigure[\label{fig:res3comp} resonance 3 in Fig.~\ref{fig:res_alpha_variation1} ($\omega\approx 1.633 - 0.026\iu$)]{\resizebox{0.45\textwidth}{!}{\includegraphics{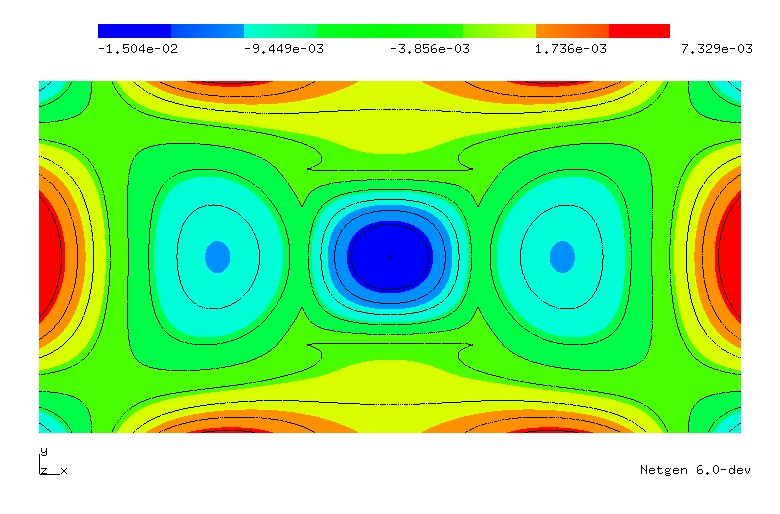}}
\resizebox{0.45\textwidth}{!}{\includegraphics{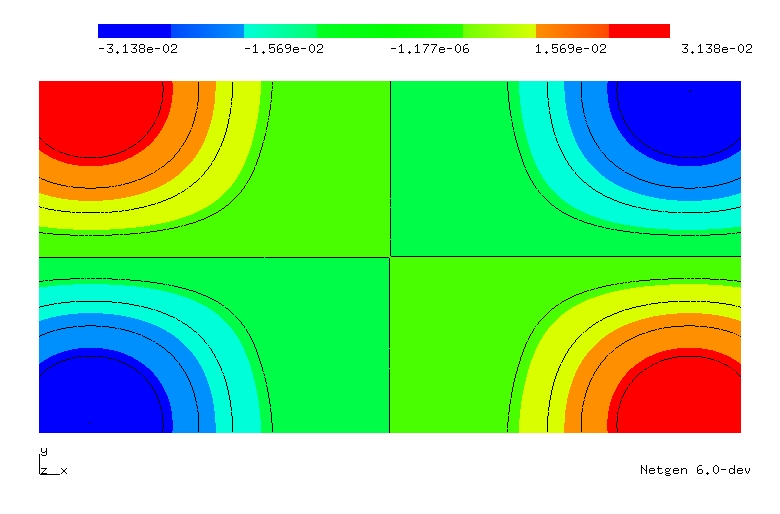}}}
\caption{real part of first (left panels) and second (right panels) Cartesian component of a resonance function for $\alpha=11$.}
\label{fig:rescomp}
\end{figure}

\dfn{In this subsection we want to highlight that the choice of the parameters $\pol,\poli$
for solving a resonance problem corresponds to a choice of the Riemann sheet on which resonances are sought.
We construct an example problem parameterized by a parameter $\alpha\in\setR^+$,
such that we have some apriori knowledge about the resonances
and observe the influence of the choice of Riemann sheets on the resonances.}

\dfo{In this subsection we }
\dfn{We} are solving \eqref{eq:res_classic} for $\Omega=\setR \times (-1,1)$
with density $\rho=4$ for $\bx\in \Omega_1:=[-0.5,0.5]^2$ and $\rho=1$ elsewhere, 
Young's modulus $E=1$ and Poisson's ratio $\nu=0.2$ and $\Bdv \bu:=\sigma(\bu)\cdot \bn$. 
Moreover, we added to \eqref{eq:res_classic} a jump in the normal stress on $\partial \Omega_1$
$$\left(\sigma(\bu^+) - \sigma(\bu^-)\right)\cdot \bn = \alpha \bu\qquad \text{for }\bx \in \partial \Omega_1,\quad \alpha \in \setC,$$
leading to an additional term $\alpha \int_{\partial \Omega_1} \bui \cdot \bvi\, ds$ \dfo{in the definition of the bilinear form $a$ in \eqref{eq:defBa}}. 

\dfn{For the discretization we have used a high order finite element method for $\Oi=(-2,2)\times (-1,1)$ based on a triangulation with maximal meshsize $h=0.1$ and a finite element order $10$. The waveguides are
discretized using the Hardy space infinite elements with $150$ basis functions in $H^-(\curve_{\pol,\poli})$ and $\pol=-0.182086-0.237784 \iu$, $\poli=-1.71492+2.29978\iu$. Since there exists a backward propagating mode (see  Fig.~\ref{fig:wn} for $\omega=1.615$), the parameters are chosen such that the pole condition is equivalent to the modal radiation condition in the neighborhood of $1.615$.}

For $\alpha \to \infty$ \dfo{this}\dfn{the additional} term leads to two decoupled problems for $\Omega_1$ and $\Omega\setminus \Omega_1$ with Dirichlet boundary conditions at $\partial \Omega_1$. Hence,
for $\alpha\to \infty$ some of the resonances $\omega_\alpha$ should converge to the square root of the (real) Dirichlet-eigenvalues of the problem in the bounded domain $\Omega_1$,
which can be computed using standard finite element methods. This can be \dfn{observed}\dfo{ seen} in Fig.~\ref{fig:res_alpha_variation1}: 
The sequences of resonances labeled with $1$ and $2$ are converging for $\alpha\to \infty$ to $\approx 1.89$, 
which is the square root of a real eigenvalue with geometric multiplicity $2$. 

\dfn{In order to distinguish resonances from discretizations of the essential spectrum, we have computed as in Sec.~\ref{sec:resprob} additionally the essential spectrum of the resonance problem, which depends on $\pol$ and $\poli$.}
For $\alpha=11$ the real part of the Cartesian components of \dfo{these}\dfn{three} resonance functions can be seen in Fig.~\ref{fig:rescomp}.
\dfo{If we use the same notation as for the Lamb waves, the resonances}
\dfn{We notice that resonances} labeled with $1$ have symmetric resonance functions (as well as these labeled with $3$)
and the resonances \dfn{labeled with} $2$ have antisymmetric resonance functions. 

\begin{figure}[tb]
\centering
\includegraphics[width=0.6\textwidth]{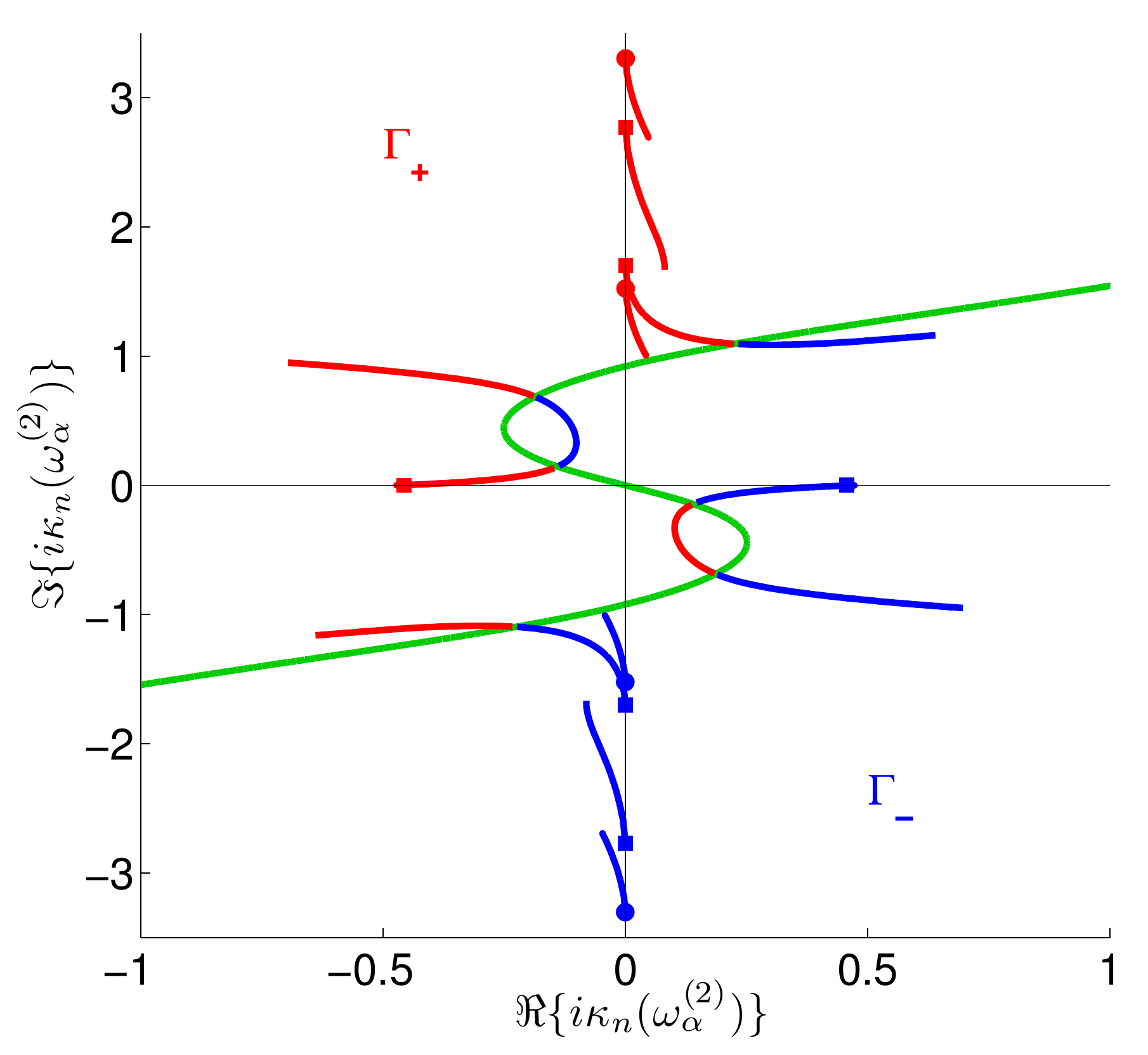}
\caption{\dfn{wavenumbers (multiplied with $\iu$) for the complex resonances 2 of Fig.~\ref{fig:res_alpha_variation1} with $\alpha \in (6,10^{10})$ computed with  \cite{Pagneux}. The squares indicate the wavenumbers for $\alpha=10^{10}$ of the symmetric Lamb waves, circles of the antisymmetric Lamb waves. The parameters for the green curve $\curve_{\pol,\poli}$ are the same as in Fig.~\ref{fig:res_alpha_variation1}.}}
\label{fig:res2_wn}
\end{figure}

\dfo{In Fig.~\ref{fig:stripholep11} the part of the spectrum with these three resonances can be seen for $\alpha=11$. There, as well as in Fig.~\ref{fig:res_alpha_variation1} there 
exists a discretization of  an essential spectrum as mentioned in Sec.~\ref{sec:elasticresprob}.}
In order to \dfn{illustrate}\dfo{study} the effect of \dfo{this}\dfn{the} essential spectrum
\dfn{on the resonance functions}, we have computed \dfn{first} for 
the sequence $2$ of \dfo{complex} resonances $\omega_\alpha^{(2)}$ in Fig.~\ref{fig:res2_wn} the wavenumbers $\kappa_n(\omega_\alpha^{(2)})$ multiplied with $\iu$.
\dfo{ of the symmetric and antisymmetric Rayleigh-Lamb modes}
The Hardy space method computes resonance functions with wavenumbers in $-\iu \curve_+$. Since for two sequences of symmetric wavenumbers there exists an intersection with $\curve$ and therefore 
jumps in these wavenumbers, we would expect jumps in the resonances $\omega_\alpha^{(2)}$, too. These jumps would be in  Fig.~\ref{fig:res_alpha_variation1} exactly at the intersection points
of $\omega_\alpha^{(2)}$ with the essential spectrum.
\dfn{Since the resonance functions corresponding to $\omega_\alpha^{(2)}$
are purely antisymmetric, they are orthogonal to the symmetric Lamb waves and thus not influenced by the change of symmetric wavenumbers. Therefore the resonances labeled with $2$ are crossing in Fig.~\ref{fig:res_alpha_variation1} the curve of the essential spectrum corresponding to a symmetric Lamb mode without any perturbation but stop at the curve of the essential spectrum corresponding to an antisymmetric Lamb mode. The last statement can be also seen in Fig.~\ref{fig:res2_wn}, since for $\alpha\approx 6$ one antisymmetric wavenumber hits $\curve$. 
}

\begin{figure}[tb]
\centering
\subfigure[\label{fig:res1_wn} $\alpha\in (5,10^{10})$. The squares indicate the wavenumbers for $\alpha=10^{10}$ of the symmetric Lamb waves, circles of the antisymmetric Lamb waves]{\resizebox{0.49\textwidth}{!}{\includegraphics{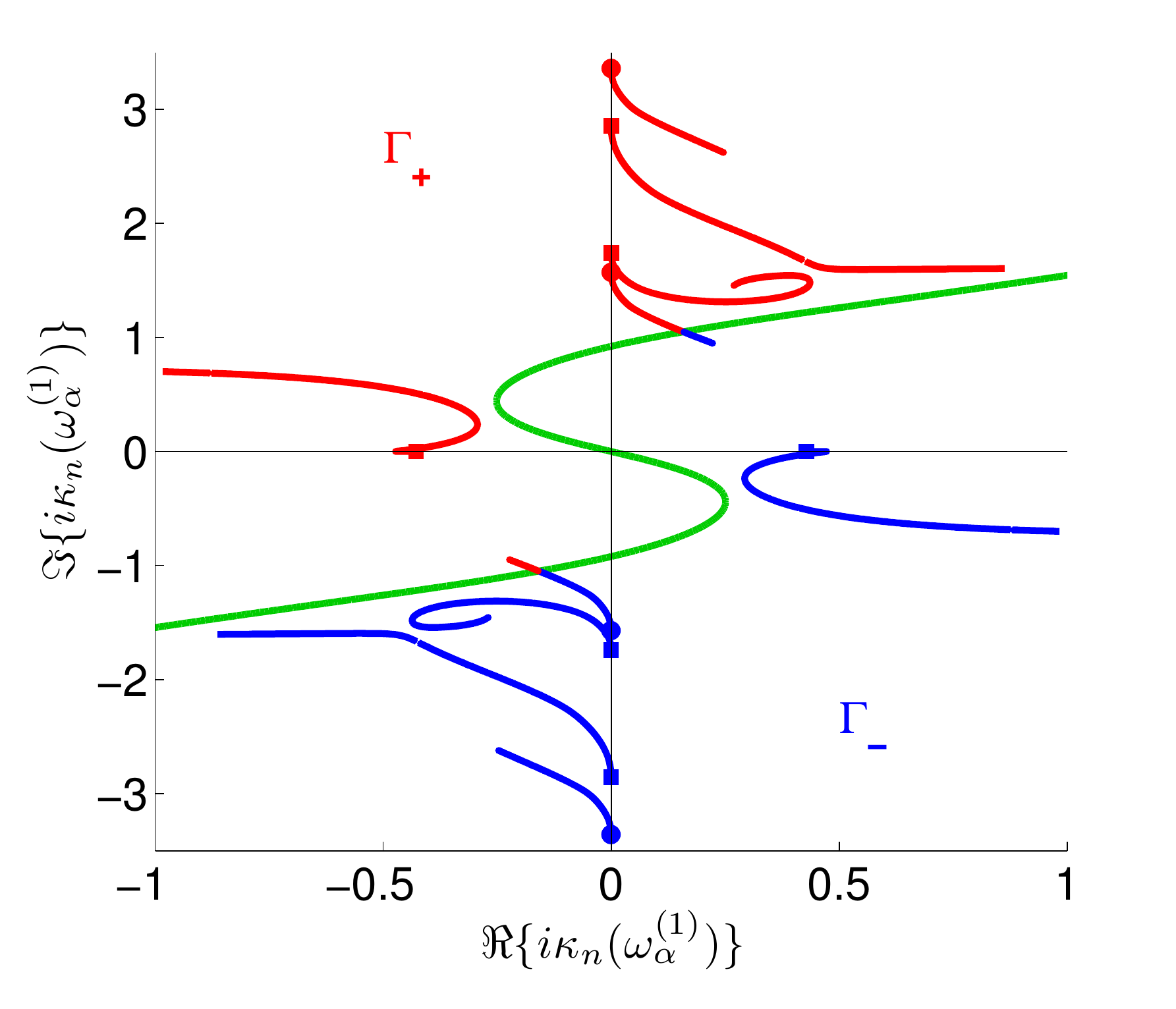}}}
\subfigure[\label{fig:res3_wn} $\alpha \in (9.4,11.7)$]{\resizebox{0.49\textwidth}{!}{\includegraphics{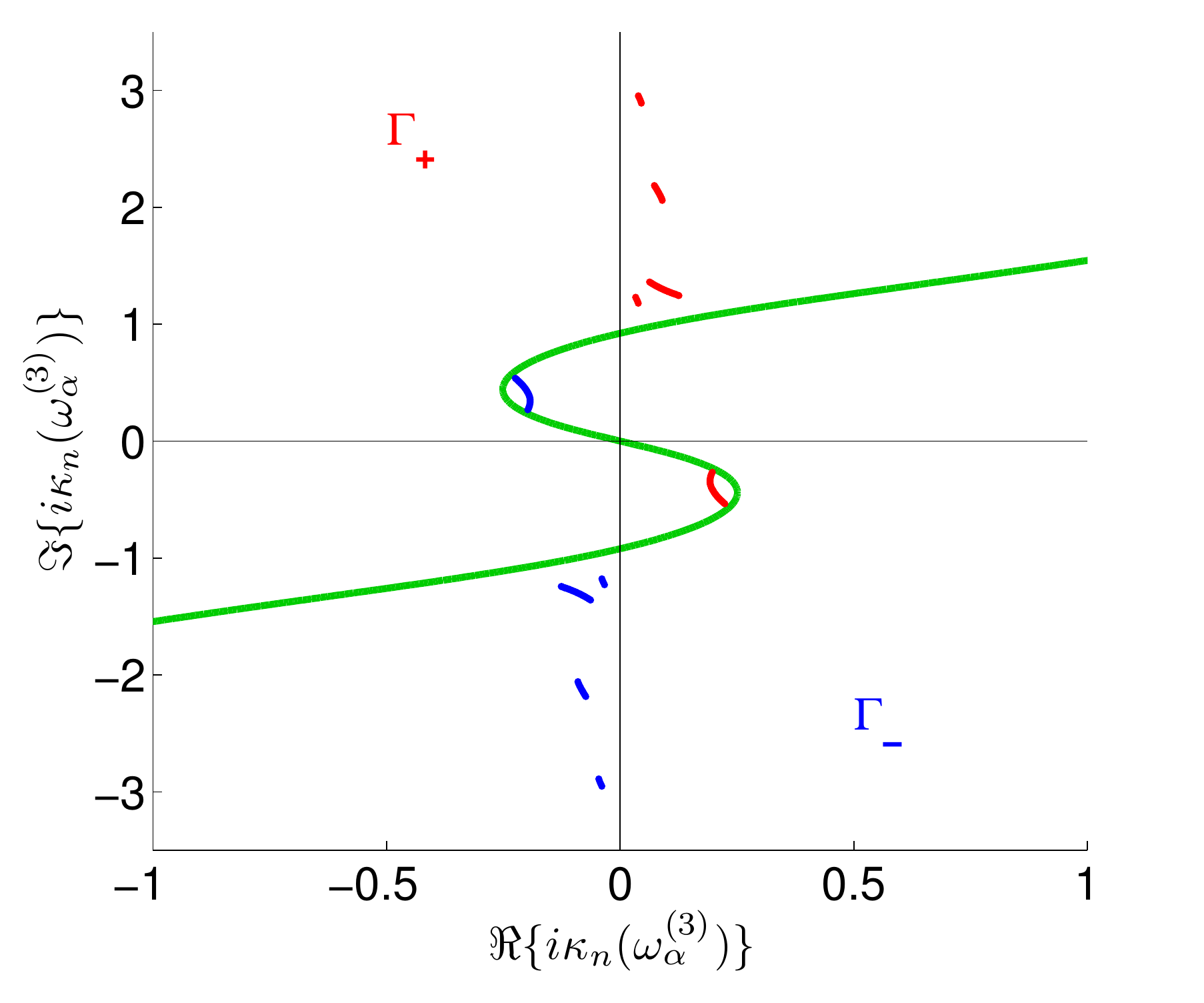}}}
\caption{wavenumbers \dfn{(multiplied with $\iu$)} for the complex resonances \dfn{1}\dfo{2} (left) and 3 (right) of Fig.~\ref{fig:res_alpha_variation1} w.r.t. $\alpha$ computed with  \cite{Pagneux}. \dfn{The parameters for the green curve $\curve_{\pol,\poli}$ are the same as in Fig.~\ref{fig:res_alpha_variation1}}.}
\label{fig:res_wn}
\end{figure}

\dfn{Fig.~\ref{fig:res_wn} shows similar results for the resonances labeled with $1$ and $3$, where the resonance functions are both symmetric. E.g. the resonances $\omega_\alpha^{(3)}$ seem to vanish after reaching the essential spectrum corresponding to a symmetric Lamb mode. The resonances $1$ are the counterparts of the resonances $3$ on a different Riemann sheet: For $\alpha=11$ both resonances are present and belong to the same Riemann sheet, which is determined by the chosen parameters $\pol$ and $\poli$. The modal radiation condition of Sec.~\ref{Sec:ModalAnalysis} for the real parts $\Re\{\omega_{11}^{(1)}\}$ and $\Re\{\omega_{11}^{(3)}\}$ of the resonances support a backward propagating mode (for the real parts of the resonances see Fig.~\ref{fig:res_alpha_variation1}, the corresponding dispersion curves are given in Fig.~\ref{fig:disprel_martin2}). $\omega_{11}^{(3)}$ is located above and $\omega_{11}^{(1)}$ below the essential spectrum arising from a branch cut of the  holomorphic extension of this backward propagating mode. For $\alpha\approx 11.7$ the resonances $\omega_{\alpha}^{(3)}$ hit the branch cut and vanish with increasing $\alpha$, since the Riemann sheet of their existence is hidden in this part of the complex plane.
Hence, for larger values of $\alpha$ we can see in Fig.~\ref{fig:res_alpha_variation1} only the resonances $\omega_{\alpha}^{(1)}$.}

\dfo{Let us now compare the resonances labeled with $1$ and $3$ in Fig.~\ref{fig:res_alpha_variation1}. The resonance functions in Fig.~\ref{fig:rescomp} look quite similar. 
Sec.~\dfo{\ref{sec:elasticresprob}}\dfn{\ref{sec:resprob}} indicates the reason for this: They are similar resonances on different \dfo{branches}\dfn{Riemann sheets} of the complex resonance problem. We expect the resonance 
sequence $\omega_\alpha^{(3)}$ to continue to the real resonance $\approx 1.89$ for $\alpha \to \infty$, if it would not be stopped by the essential spectrum, which is given by 
the branch cuts. This can be also seen in Fig.~\ref{fig:res3_wn}. Our choice of $\pol$ and $\poli$ and therefore our choice of $\curve$ determines, on which part of the complex plane 
we see which branch of the resonance problem. The resonances $\omega_\alpha^{(2)}$ are the same in two of these branches, since they are independent of the symmetric Lamb mode ``producing'' the branch.}


\dfo{In order to clarify the relevance of one resonance for a diffraction problem with positive frequency $\omega$, we have computed in Fig.~\ref{fig:striphole2} the resonances additionally for a standard 
Hardy space method, where $\curve$ is a straight line. In other words we see the resonances on other branches in comparison to $\pol=-0.182086-0.237784\iu$, $\poli=-1.71492+2.29978\iu$. 
Eigenvalues with positive imaginary part are an indicator, that something went wrong (in some sense). Indeed, for $\omega\in (1.58,1.65)$ we have modes with different signs of group and phase velocity
(case two in Sec.~\ref{Sec:CondCurve}). So, for $\omega\in (1.58,1.65)$ the choice $\pol=\poli=-1+\iu$ leads to a physically incorrect radiation condition. 

In Fig.~\ref{fig:stripholep11} we haven chosen the correct radiation condition for $\omega\in (1.58,1.65)$. Since the resonances $2$ and $3$ lie on the same branch as $\omega$, 
the relevance of these resonances can be measured by the Euclidean distance between $\omega$ and them. The resonance $3$ is below the essential spectrum and therefore belongs 
to another branch. The distance of such a resonance to $\omega$ would be given by the length of a path around the branch cut.

Note, that by the choice of parameters $\pol=-0.182086-0.237784\iu$ and $\poli=-1.71492+2.29978\iu$ the radiation condition is correct for $\omega\in (1.58,1.65)$, 
but not for all $\omega>0$. Fig.~\ref{fig:res_alpha_variation1} indicates by the existence of eigenvalues with positive imaginary part, that the radiation
condition is incorrect for $\omega<1.42$. Hence, the resonances labeled with $5$ lie on another branch as the real frequencies $\omega<1.42$. Therefore, their
relevance is also limited (see Fig.~\ref{fig:ess_spec} for a larger frequency interval).}

\section{\dfn{Conclusion}}
\dfn{In this paper we addressed the topic of backward propagating modes in time-harmonic two-dimensional elastic waveguides. The so called pole condition allows to reformulate the radiation condition even in the presence of backward waves without using the waveguide modes or wavenumbers. For the problems under consideration a special class of pole conditions depending on two complex parameters $\pol$ and $\poli$
is sufficient. Detailed explanations on the choice of these parameters are given in Section~\ref{Sec:CondCurve}. The Hardy space infinite element method developed in \cite{HallaHohageNannenSchoeberl:14} relies on this pole condition and can be easily applied to waveguide structures using tensor product elements.}

\dfn{An outstanding property of the presented scheme is that it is independent of the waveguide modes and wavenumbers which simplifies the implementation. Moreover, a discetization of a resonance problem leads to a generalized linear matrix eigenvalue problem, which is much easier to handle than a non-linear problem. However, the interpretation of the occurring spectral objects
for resonance problems is rather involved and was addressed in Section~\ref{sec:resprob} and Section~\ref{sec:MP3}.}

\dfn{The method shows in numerical examples super-algebraic convergence for diffraction problems. For resonance problems the numerical results also coincide with the theoretical considerations. In particular, a strong relation between the resonances and the behavior of the solutions to diffraction problems with frequencies in the neighborhood of a resonance frequency can be seen. Hence, resonance and diffraction problems in two dimensional waveguides can be reliably solved using the Hardy space infinite element method.}

\bibliographystyle{siam}
\bibliography{bibliography}

\begin{thebibliography}{10}

\bibitem{Abarbanel:99}
{\sc S.~Abarbanel, D.~Gottlieb, and J.~Hesthaven}, {\em Well-posed perfectly
  matched layers for advective acoustics}, Journal of Computational Physics,
  154 (1999), pp.~266 -- 283.

\bibitem{Achenbach73}
{\sc J.~D. Achenbach}, {\em {Wave Propagation in Elastic Solids (North-Holland
  Series in Applied Mathematics and Mechanics)}}, North-Holland series in
  applied mathematics and mechanics, v. 16, North Holland, 1987.

\bibitem{Baronian}
{\sc V.~Baronian, A.-S. Bonnet-BenDhia, and E.~Lun{\'e}ville}, {\em
  {Transparent boundary conditions for the harmonic diffraction problem in an
  elastic waveguide}}, Journal of Computational and Applied Mathematics, 234
  (2010), pp.~1945--1952.

\bibitem{BNiBonLeg:2004}
{\sc {\'E}.~B{\'e}cache, A.-S. Bonnet-BenDhia, and G.~Legendre}, {\em Perfectly
  matched layers for the convected {H}elmholtz equation}, SIAM Journal on
  Numerical Analysis, 42 (2004), pp.~409--433.

\bibitem{Joly:03}
{\sc E.~B{\'e}cache, S.~Fauqueux, and P.~Joly}, {\em Stability of perfectly
  matched layers, group velocities and anisotropic waves}, J. Comput. Phys.,
  188 (2003), pp.~399--433.

\bibitem{Berenger:94}
{\sc J.-P. Berenger}, {\em A perfectly matched layer for the absorption of
  electromagnetic waves}, J. Comput. Phys., 114 (1994), pp.~185--200.

\bibitem{NEVP}
{\sc W.-J. Beyn}, {\em An integral method for solving nonlinear eigenvalue
  problems}, Linear Algebra and its Applications, 436 (2012), pp.~3839 -- 3863.
\newblock Special Issue dedicated to Heinrich Voss's 65th birthday.

\bibitem{Bonnet}
{\sc A.-S. Bonnet-BenDhia, C.~Chambeyron, and G.~Legendre}, {\em On the use of
  perfectly matched layers in the presence of long or backward guided elastic
  waves}, Wave Motion,  (to appear).

\bibitem{Braess:FE}
{\sc D.~Braess}, {\em Finite Elemente}, Springer, 2003.

\bibitem{Chimenti:97}
{\sc D.~E. Chimenti}, {\em Guided waves in plates and their use in materials
  characterization}, Applied Mechanics Reviews, 50 (1997), pp.~247--284.

\bibitem{Givoli:04}
{\sc D.~Givoli}, {\em High-order local non-reflecting boundary conditions: a
  review}, Wave Motion, 39 (2004), pp.~319--326.

\bibitem{Graff}
{\sc K.~Graff}, {\em Wave Motion in Elastic Solids}, Oxford engineering science
  series, Clarendon Press, 1975.

\bibitem{HallaHohageNannenSchoeberl:14}
{\sc M.~Halla, T.~Hohage, L.~Nannen, and J.~Sch\"oberl}, {\em Hardy space
  infinite elements for time-harmonic wave equations with phase velocities of
  different signs}, preprint, Institute for Analysis and Scientific Computing,
  TU Wien, 2014.

\bibitem{HohageNannen:09}
{\sc T.~Hohage and L.~Nannen}, {\em Hardy space infinite elements for
  scattering and resonance problems}, SIAM J. Numer. Anal., 47 (2009),
  pp.~972--996.

\bibitem{HohageNannen:15}
\leavevmode\vrule height 2pt depth -1.6pt width 23pt, {\em Convergence of
  infinite element methods for scalar waveguide problems}, BIT Numerical
  Mathematics, 55 (2015), pp.~215--254.

\bibitem{PC1}
{\sc T.~Hohage, F.~Schmidt, and L.~Zschiedrich}, {\em Solving time-harmonic
  scattering problems based on the pole condition. {I}. {T}heory}, SIAM J.
  Math. Anal., 35 (2003), pp.~183--210.

\bibitem{Orazov}
{\sc A.~G. Kostyuchenko and M.~B. Orazov}, {\em Problem of oscillations of an
  elastic half cylinder and related self-adjoint quadratic pencils}, Journal of
  Soviet Mathematics, 33 (1986), pp.~1025--1065.

\bibitem{Lamb:1904}
{\sc H.~Lamb}, {\em On group - velocity}, Proceedings of the London
  Mathematical Society, s2-1 (1904), pp.~473--479.

\bibitem{Langenberg:12}
{\sc K.-J. Langenberg, R.~Marklein, and K.~Mayer}, {\em Ultrasonic
  Nondestructive Testing of Materials}, {CRC} Press, feb 2012.

\bibitem{LM:08}
{\sc M.~Levitin and M.~Marletta}, {\em A simple method of calculating
  eigenvalues and resonances in domains with infinite regular ends}, Proc. Roy.
  Soc. Edinburgh Sect. A, 138 (2008), pp.~1043--1065.

\bibitem{Lighthill:65}
{\sc M.~J. Lighthill}, {\em Group velocity}, J. Inst. Math. Appl., 1 (1965),
  pp.~1--28.

\bibitem{SMARTlayer2}
{\sc L.~Métivier, R.~Brossier, S.~Labbé, S.~Operto, and J.~Virieux}, {\em A
  robust absorbing layer method for anisotropic seismic wave modeling}, Journal
  of Computational Physics, 279 (2014), pp.~218 -- 240.

\bibitem{Nannenetal:13}
{\sc L.~Nannen, T.~Hohage, A.~Sch{\"a}dle, and J.~Sch{\"o}berl}, {\em Exact
  {S}equences of {H}igh {O}rder {H}ardy {S}pace {I}nfinite {E}lements for
  {E}xterior {M}axwell {P}roblems}, SIAM J. Sci. Comput., 35 (2013),
  pp.~A1024--A1048.

\bibitem{ngs-waves}
{\sc L.~Nannen and J.~Sch{\"o}berl}, {\em Software module ngs-waves}.
\newblock http://sourceforge.net/projects/ngs-waves/, 2014.
\newblock addon to the mesh generator Netgen and the high order finite element
  code NGSolve.

\bibitem{Pagneux}
{\sc V.~Pagneux and A.~Maurel}, {\em Determination of lamb mode eigenvalues},
  The Journal of the Acoustical Society of America, 110 (2001), pp.~1307--1314.

\bibitem{netgen}
{\sc J.~Sch\"oberl}, {\em Netgen - an advancing front 2d/3d-mesh generator
  based on abstract rules}, Comput.Visual.Sci, 1 (1997), pp.~41--52.

\bibitem{ngsolve:14}
\leavevmode\vrule height 2pt depth -1.6pt width 23pt, {\em C++11 implementation
  of finite elements in ngsolve}, Preprint 30/2014, Institute for Analysis and
  Scientific Computing, TU Wien, 2014.

\bibitem{Skeltonetal:07}
{\sc E.~A. Skelton, S.~D.~M. Adams, and R.~V. Craster}, {\em Guided elastic
  waves and perfectly matched layers}, Wave Motion, 44 (2007), pp.~573--592.

\bibitem{SMARTlayer1}
{\sc J.~Tago, L.~Métivier, and J.~Virieux}, {\em Smart layers: a simple and
  robust alternative to pml approaches for elastodynamics}, Geophysical Journal
  International, 199 (2014), pp.~700--706.

\end{thebibliography}
\end{document}